\documentclass[a4paper,10pt]{amsart}
\usepackage[english]{babel}
\usepackage{amsmath,amsthm,amssymb,amsfonts,graphicx,esint}

\usepackage{enumitem}
\setlist[enumerate]{leftmargin=*}
\setlist[itemize]{leftmargin=*}

\usepackage{mathrsfs}

\usepackage{hyperref}
\hypersetup{
  colorlinks = true,
  urlcolor   = magenta,
  linkcolor  = blue,
  citecolor  = red
}

\numberwithin{equation}{section}

\newtheorem{theorem}{Theorem}[section]
\newtheorem{lemma}[theorem]{Lemma}
\newtheorem{corollary}[theorem]{Corollary}
\newtheorem{proposition}[theorem]{Proposition}
\theoremstyle{definition}
\newtheorem{definition}[theorem]{Definition}

\newtheorem{remark}[theorem]{Remark}

\DeclareMathOperator*{\esssup}{ess\,sup}
\newcommand{\supp}{\mathop{\mathrm{supp}}}
\newcommand{\dist}{\mathrm{dist}}
\newcommand{\tonde}[1]{\left( #1 \right)}   
\newcommand{\quadre}[1]{\left[ #1 \right]}   
\newcommand{\graffe}[1]{\left\{ #1 \right\}} 
\newcommand{\abs}[1]{\left| #1 \right|}      
\newcommand{\norma}[1]{\left\| #1 \right\|}  

\newcommand{\argomento}{\vcenter{\hbox{\scalebox{1.5}{$\cdot$}}}}
\newcommand{\RR}{\mathbb{R}}
\newcommand{\Rz}{\mathcal{R}}

\newcommand{\ZZ}{\mathbb{Z}}
\newcommand{\NN}{\mathbb{N}}

\newcommand{\Q}{\mathcal{Q}} 
\newcommand{\vfX}{X}
\newcommand{\vfXN}{X^\flat}
\renewcommand{\P}{\mathcal{P}}

\newcommand{\ee}{\mathrm{e}}         
\newcommand{\dd}{\mathrm{d}}          
\newcommand{\lie}[1]{\mathfrak{#1}}   
\newcommand{\dil}{\mathfrak{d}}       
\newcommand{\dimV}{d_1}               
\newcommand{\dimZ}{d_2}               
\newcommand{\dimHomog}{Q}
\newcommand{\subriemdist}{\mathrm{d}_{\mathrm{SR}}} 
\newcommand{\subriemball}[2]{\mathrm{B}_{\mathrm{SR}}\tonde{{#1},{#2}}}
\newcommand{\riemdist}{\mathrm{d}_{\mathrm{R}}} 
\newcommand{\riemball}[2]{\mathrm{B}_{\mathrm{R}}\tonde{{#1},{#2}}}
\newcommand{\damekriccidist}{\mathrm{d}_{\mathrm{S}}}
\newcommand{\damekricciball}[2]{\mathrm{B}_{\mathrm{S}}\tonde{{#1},{#2}}}
\newcommand{\diltile}[2]{\dil_{4^{#1}}\tonde{#2}}
\newcommand{\axbdist}{\mathrm{d}_{ax+b}}
\newcommand{\horiball}[2]{\mathrm{B}_{\mathfrak{v}}\tonde{{#1},{#2}}}
\newcommand{\vertball}[2]{\mathrm{B}_{\mathfrak{z}}\tonde{{#1},{#2}}}
\newcommand{\diam}{\mathrm{diam}}

\newcommand{\supeuclball}[2]{\mathrm{B}_{\infty}\tonde{{#1},{#2}}}
\newcommand{\costante}{{\tau_N}}
\newcommand{\costantee}{\mathcal K_\ast }
\newcommand{\constantS}{\tau}
\newcommand{\Nprod}{\cdot_{N}}
\newcommand{\rI}{r_I}
\newcommand{\Sprod}{\cdot_{S}}
\newcommand{\modular}{\delta}
\newcommand{\riesz}{\mathcal R}
\newcommand{\children}{\mathcal{C}}
\newcommand{\admis}{\mathcal A}

\newcommand{\indicator}[1]{\mathbf{1}_{#1}}
\newcommand{\BMO}{\mathrm{BMO}}

\newcommand{\bx}{\mathbf{x}}
\newcommand{\by}{\mathbf{y}}
\newcommand{\bz}{\mathbf{z}}

\newcommand{\tc}{\mathrel{\,:\,}}     
\newcommand{\defeq}{\mathrel{:=}}     
\newcommand{\eqdef}{\mathrel{=:}}

\title[Calder\'on--Zygmund Theory on Damek--Ricci spaces]{Calder\'on--Zygmund Theory and Endpoint Results for Riesz Transforms on Damek--Ricci Spaces}
\author{Martina Bernardelli}
\address[M. Bernardelli]{Dipartimento di Scienze Matematiche ``Giuseppe Luigi Lagrange'' \\ Politecnico di Torino \\ Corso Duca degli Abruzzi 24 \\ 10129 Torino \\ Italy}
\email{martina.bernardelli@polito.it}

\author{Jie Liu}
\address[J. Liu]{School of Mathematics and Statistics \\ Northwestern Polytechnical University \\ Xi'an 710129 \\ China \\ and Dipartimento di Scienze Matematiche ``Giuseppe Luigi Lagrange''\\ Politecnico di Torino \\ Corso Duca degli Abruzzi 24 \\ 10129 Torino \\ Italy}
\email{jay2000@mail.nwpu.edu.cn}
\email{jie.liu@polito.it}


	
\thanks{The first-named author is member of Gruppo Nazionale per l'Analisi Matematica, la Probabilit\`a e le loro Applicazioni (GNAMPA) of Istituto Nazionale di Alta Matematica (INdAM). The second-named author gratefully acknowledges the financial support of the China Scholarship Council (Grant No.\ 202406290140).}

\begin{document}
\begin{abstract}
We develop an atomic Hardy space and Calder\'on--Zygmund theory on Damek--Ricci spaces equipped with a left-invariant Riemannian distance and a right Haar measure. The main geometric ingredient is the construction of a regular Calder\'on--Zygmund structure adapted to the two-scale geometry of these spaces. Starting from Strichartz's self-similar tiling of the underlying $H$-type group $N$, we construct an explicit system of cubes on $N$ compatible with the Riemannian metric. Combining these cubes with dyadic intervals in $\RR$ yields a suitable family of admissible sets in the Damek--Ricci space $S=N\rtimes \RR$ satisfying the required recursive decomposition, exhaustion, regularity, and enlargement properties. This provides a concrete atomic Hardy space $H^1(S)$ and a bounded mean oscillation space $\mathrm{BMO}(S)$ with the natural duality and interpolation properties. We then study the first-order Riesz transforms associated with the distinguished Laplacian. We prove that such Riesz transforms are bounded from \(H^1(S)\) to \(L^1(S)\), whereas their adjoint fail to have this endpoint boundedness. By \(H^1\)--BMO duality, it follows that the Riesz transforms are not bounded from \(L^\infty(S)\) to \(\mathrm{BMO}(S)\) while their adjoint are bounded. Thus the Riesz transforms and their adjoints, although bounded on \(L^p(S)\) for \(1<p<\infty\), exhibit an
asymmetric behaviour at the endpoints.
\end{abstract}

\maketitle

\section{Introduction}
Calder\'on--Zygmund (CZ) theory, originating in the work of Calder\'on and Zygmund in the 1950s, provides a fundamental framework for the study of singular integral operators. In the classical Euclidean setting, its mechanism is based on the existence of the so-called CZ decomposition of integrable functions. This decomposition interacts well with the kernel of CZ operators, thanks to some properties of $\RR^n$, such as the doubling property and the existence of a dyadic structure. Combined with localization, cancellation and off-diagonal kernel regularity, it leads to weak type $(1,1)$ and $L^p$-boundedness, as well as endpoint estimates on Hardy and $\BMO$ spaces; see, for instance, \cite{CZ52,C74}.

It was soon clear to the community of mathematicians that the rigid structure of $\RR^n$ was not necessarily needed to develop a CZ theory. Indeed, Coifman and Weiss in the 1970s extended the theory to spaces of homogeneous type, where the doubling condition provides the measure control required by the classical arguments \cite{CW71,CW77}. At that time a system of cubes analogous to the dyadic structure in $\RR^n$ was not available, but they figured out that the theory was flexible enough to be applied to balls instead of cubes. Systems of dyadic cubes adapted to general doubling metric spaces were later constructed by Christ in 1990 \cite{C90} and, in a more flexible form, by Hyt\"onen and Kairema in 2012 \cite{HK12}.

Beyond spaces of homogeneous type, CZ theory has also been developed in various nondoubling settings. In particular, in the early 2000s Tolsa and, independently and in parallel, Nazarov, Treil and Volberg extended the theory to $\RR^n$ endowed with a non necessarily locally doubling measure, satisfying polynomial growth conditions. In this setting the endpoint estimates are expressed through the atomic Hardy space and the space $\mathrm{BMO}$ introduced by Tolsa, which play the role of the classical endpoint spaces \cite{NTV02,T01}.

Another line of research dealt with the global failure of the doubling property: large metric balls do not satisfy a uniform doubling condition, and therefore cannot serve as localization sets in the usual CZ arguments. This kind of setting has attracted considerable interest over the last 30 years; see, for example \cite{GS96}. One of the most important contributions in this direction is the work of Hebisch and Steger \cite{HS03} in 2003, who developed a unified CZ theory on some Lie groups with exponential growth based on suitably adapted admissible sets. One of the novelties of their work lies in the introduction of a system of admissible sets, not necessarily balls or cubes, but comparable to them in a proper sense, that compensated for the lack of the doubling property. Their approach was subsequently used to construct atomic Hardy and $\mathrm{BMO}$ spaces on the $ax+b$ group and on more general solvable extensions of stratified groups; see, for example, \cite{LVY12,MOV18,V07,V09}.

A very recent development of the theory is due to Conde Alonso, Martini and Vallarino \cite{CMV26}, who generalized the framework to abstract measure spaces, not necessarily endowed with a metric, provided that they support a family of admissible sets satisfying properties such as nesting (the existence of an at most countable family of children that partition each admissible set), exhaustion (the existence of an increasing sequence of admissible sets which cover the whole space), and a suitable approximation of the measurable sets by admissible sets.

Our first aim is to develop a CZ and Hardy space theory on Damek--Ricci spaces. To this end, we construct a system of admissible sets, in the spirit of the general theory developed in \cite{CMV26}.

Damek--Ricci spaces were introduced as harmonic extensions of $H$-type groups by Damek and Ricci \cite{D87a,D87b,DR92a,DR92b} and studied intensely thereafter; see, for example, \cite{ADY96,CDKR91,CDKR98}. To define such spaces we need to introduce $H$-type algebras. An $H$-type algebra in the sense of Kaplan \cite{K90} is a $2$-step nilpotent Lie algebra $\lie n=\lie v\oplus\lie z$ equipped with an inner product such that, if $\lie v$ is the orthogonal complement of the centre $\lie z$ of $\lie n$, and $\lie z^*$ is the dual of $\lie z$, then
\begin{equation*}
\abs{\theta[x,\cdot]}=\abs{\theta}\abs{x}\qquad\forall x\in \lie v\ \ \forall\mu\in \lie z^*,
\end{equation*}
where $\theta[x,\cdot]$ is thought of as an element of the dual $\lie v^*$ of $\lie v$. The direct sum decomposition $\lie n=\lie v\oplus\lie z$ is a stratification of $\lie n$, since $[\lie v,\lie v]=\lie z$ (when $\lie v$ is nontrivial). We shall refer to $\lie v$ as the \emph{horizontal layer} and to $\lie z$ as the \emph{vertical layer}, and we shall use the same terminology for the corresponding layers of $N$, the connected, simply connected Lie group associated with $\lie n$. We denote by $\dimV$ and $\dimZ$ the dimensions of the layers $\lie v$ and $\lie z$, and by $\graffe{\vfXN_{j}}_{j=1}^{\dimV+\dimZ}$ the frame of left-invariant vector fields on $N$ such that $\{\vfXN_{j}(0,0)\}_j$ is an orthonormal basis of $\lie n$ adapted to the decomposition, that is, $\vfXN_{j}(0,0)\in \lie v$ for $j=1,\dots,\dimV$ and $\vfXN_{j}(0,0)\in \lie z$ for $j=\dimV+1,\dots,\dimV+\dimZ$ where $(0,0)$ is the identity element of $N$. We equip $N$ with the natural family $\graffe{\dil_t\tonde{\cdot}}_{t\in\RR^+}$ of automorphic dilations acting on $\lie n$ by
\begin{equation}\label{Def-dilation}
\dil_t{(x,z)}=\big(t^{\frac12}x,\,tz\big)\qquad\forall (x,z)\in \lie v\oplus \lie z,\,\,\,\forall t\in\RR^+.
\end{equation} 
More details about $H$-type groups can be found in Subsection \ref{subsection-Damek–Ricci spaces}.

A Damek--Ricci space $S$ is the semidirect extension $S=N\rtimes\RR$, obtained by letting the additive group $\RR$ act on $N$ through the automorphic dilations $\dil_{\ee^u}$, via the exponential coordinate $u$ on the extension factor.

The Lie algebra $\lie s$ of $S$ is then the semidirect extension $\lie n\rtimes\RR$, where $\partial_u$ acts on $\lie n$ as the derivation generating the dilations $\dil_{\ee^u}$. Let $\vfX_0$ be the generator of the $\RR$ direction, then $\graffe{\vfX_j(e)}_{j=0}^{\dimV+\dimZ}$ is an orthonormal basis of $\lie s$, with $\vfX_j(e)=\vfXN_j(0,0)$ for $j\geq1$ where $e$ is the identity element of $S$. Transporting the basis around $S$ by left translations, we obtain a frame of left-invariant vector fields such that
\begin{equation} \label{Def-fields}
\begin{aligned}
    \vfX_0&=\partial_u,\\
    \vfX_j&=\ee^{u/2}\vfXN_j
    && j=1,\dots,\dimV,\\
    \vfX_{\dimV+j}&=\ee^u\vfXN_{\dimV+j}
    && j=1,\dots,\dimZ.
\end{aligned}
\end{equation}
More details about those fields will appear in Section \ref{Sec-endpointRiesz}.

We endow $S$ with the left-invariant Riemannian metric for which this frame is orthonormal at every point. Let $\lambda$ and $\rho$ denote left and right Haar measures on $S$, which are different since $S$ is not unimodular. Both the $\rho$-measure and the $\lambda$-measure of geodesic balls grow exponentially in the radius, making $S$ a group of exponential volume growth. Following the literature \cite{CGHM94, GS96, HS03, H74, LM26, MV10, V07} among others, we equip $S$ with the left-invariant Riemannian metric $\damekriccidist$ and the right Haar measure $\rho$, and hence for all $p\in[1,\infty]$, we denote by $L^p(S)$ the $L^p$ space with respect to $\rho$.

For every $u\in\RR$, the slice $N_u\defeq N\times\graffe{u}$ is a submanifold of $S$ whose tangent space at $(n,u)\in N_{u}$ is spanned by $\graffe{\vfX_j(n,u)}_{j=1}^{\dimV+\dimZ}$, which forms a left-invariant orthonormal frame for the  Riemannian metric tensor of $S$ restricted to $N_u$; 
we denote its intrinsic distance by $\dd_{R,u}$. For $u=0$, we get the Riemannian metric of $N$ for which $\graffe{\vfXN_j}$ is orthonormal, and for general $u$ we obtain an expression of the distance from it by dilation, namely
\[
\dd_{R,u}(n_1,n_2)=\dd_{R,0}\big(\dil_{\ee^{-u}}n_1,\dil_{\ee^{-u}}n_2\big), \qquad n_1,n_2\in N.
\]
The Riemannian distance $\dd_{R,0}$ on $N$ is equivalent at small scales to the Euclidean distance, and at large scales to the intrinsic Carnot--Carath\'eodory left-invariant (sub-Riemannian) distance $\subriemdist$ associated with the horizontal vector fields $\{\vfXN_j\}_{j=1}^{\dimV}$. Also $\dd_{R,u}$ exhibits the same two regimes, but the scale at which one passes from one to the other depends on $u$.

The sub-Riemannian distance $\subriemdist$ on $N$ is left-invariant and homogeneous with respect to the group dilations, that is,
\[
\subriemdist(\dil_t n_1,\dil_t n_2)=t^{1/2}\,\subriemdist(n_1,n_2) \,\,\,\text{ for all }n_1,n_2\in N.
\]
The Riemannian distance $\riemdist$, however, does not enjoy this homogeneity. Therefore, although the Riemannian geometry is the one relevant to our purposes, the sub-Riemannian distance provides a useful tool at large scales, thanks to its compatibility with the group dilations. The details will appear in the construction of the cubes on $N$ in Section \ref{Sec-dyadicinN}.

This dichotomy between large and small scales appears also on $S$, and it is already present in the earlier literature on Damek--Ricci spaces: for the Riemannian metric, it has proved difficult to find a single system of admissible sets that works at all scales. Following the approach of Hebisch and Steger in \cite{HS03}, Vallarino \cite{V07} developed a CZ theory on Damek--Ricci spaces and applied it to establish weak-type $(1,1)$ estimates for spectral multipliers of the distinguished Laplacian. Her construction relies on admissible sets of two different types: at large scales products of dyadic sets in $N$ and intervals in $\RR$, and at small scales geodesic balls on $S$.

Our construction retains the product structure at both large and small scales, differing substantially to Vallarino's one at small scales. Indeed, the dyadic cubes in $N$, constructed in \cite{V07} using Christ's theorem \cite[Theorem 11]{C90} as a black box, are adapted to the sub-Riemannian distance and therefore do not directly provide the uniform geometric control needed at arbitrarily small Riemannian scales. In particular, the comparison between these cubes and their metric enlargements becomes problematic. This difficulty is avoided in \cite{V07} by using geodesic balls at small scales, which are sufficient for the CZ theory developed there. However, such a family does not provide the nested dyadic partitions required by the framework of \cite{CMV26}, which is suited to establishing Hardy-space boundedness for Riesz transforms and other CZ operators satisfying the appropriate hypotheses.

Another class of admissible sets on Damek--Ricci spaces (and more general rank-one extensions of Carnot groups) was introduced by Martini, Ottazzi and Vallarino \cite{MOV18}, who took as admissible sets the products of Christ cubes on $N$ with intervals of $\RR$. Their construction, though, is adapted to the sub-Riemannian geometry of the group, since they worked with the sub-Laplacian on $S$, and is therefore not suitable for our purposes, at least at small scales, where the sub-Riemannian geometry of $S$ differs from the Riemannian one.

To obtain a nested system of sets in $S$ compatible with the Riemannian metric, we introduce a new family of Riemannian dyadic cubes at small scales on $N$. These sets are constructed from sub-Riemannian cubes on $N$, retain an explicit nested subdivision structure, and are uniformly comparable to Riemannian $N$-balls. 

We are able to carry out the above structure thanks to Strichartz, who in 1992 constructed an explicit self-similar tiling of $H$-type groups; see \cite[Section 2]{S92}. The explicit algebraic form he provides yields a system of cube-like sets on $N$, and allows us to show directly how (small) sub-Riemannian cubes can be ``glued'' into Riemannian ones, so that the passage between the two geometries preserves the nesting properties and the approximation of measurable sets we wanted, while keeping under control the comparison between the measures of cubes and their metric enlargements. Moreover, this passage can be performed at every scale, which is necessary because on $S$ the $\RR$ coordinate determines the geometry of the $N$-component by dilation.

The admissible sets on \(S=N\rtimes\RR\) are therefore products
\[
R=\Q\times I, 
\]
where \(\Q\) is a dyadic cube on \(N\) and $I$ is an interval of $\RR$. The $N$ and $\RR$ scales cannot be chosen independently: the \(\RR\)-variable acts on the \(N\)-component by dilation, so that the effective size of \(\Q\) depends on the midpoint of \(I\), while the variation across \(I\) adds a further dependence on the length of $I$. We therefore impose an admissibility condition that balances these quantities (see Section \ref{Sect-adminS}). Partitioning $R$ in a suitable way into other admissible sets yields its children $\children(R)$; iterating this procedure $n$ times, we obtain the $n$-th generation of children $\children^n(R)$ of $R$. The family of all such generations is called the family of descendants of $R$ and is denoted by $\children^\ast(R)$.

To use those sets in CZ arguments, however, one also needs uniform control of their behaviour under metric enlargement. Given an admissible set \(R\), we define
\begin{equation}\label{Def-enlargedset}
R^\ast\defeq\left\{\mathbf x\in S:\damekriccidist(\mathbf x,R)\leq\constantS\cdot\diam_{\mathrm S}(R)\right\},
\end{equation}
where \(\constantS>0\) is a fixed constant which will be presented explicitly in Section \ref{Sect-adminS} and 
\begin{equation*}
\diam_{\mathrm S}(R)\defeq \sup_{x,y\in R} \mathrm{d}_{\mathrm{S}}(x,y).
\end{equation*}
The key property is that the enlargement does not increase the measure by more than a constant factor, namely, there exists $C>0$ such that
\begin{equation}\label{Eq-enlargementestimate}
\rho(R^\ast)\leq C \rho(R)\quad \text{for every admissible set } R.
\end{equation}
In the absence of global doubling, \eqref{Eq-enlargementestimate} controls the measure of the union of the enlarged admissible sets arising in the CZ decomposition and thus provides the localization needed in singular integral estimates.

Collecting the recursive parent-children structure, the exhaustion properties and the enlargement estimate \eqref{Eq-enlargementestimate}, we obtain the following result, where \(\admis\) is the family of admissible sets whose precise construction will appear in Section \ref{Sect-adminS}.

\begin{theorem}\label{Thm-CZspaceproof}
    The Damek--Ricci space $S$, endowed with the left-invariant Riemannian metric $\damekriccidist$, the right Haar measure $\rho$ on Borel sets and the system of admissible sets $\mathcal A$, is a regular metric CZ space in the sense of \cite[Definition~4.1]{CMV26}, that is:
    \begin{enumerate}[label=(\alph*)]
        \item\label{regular a} each admissible set has finite non-zero measure;
        \item\label{regular b} there exists an \emph{exhausting sequence} $(R_n)_{n\in \NN}$ of admissible sets, satisfying $R_n\subseteq R_{n+1 }$ for all $n\in \NN$ and $S=\bigcup_{n\in\NN}R_n;$
        \item\label{regular c} each admissible set $R$ partitions into a finite family $\children(R) \subseteq \admis$ whose elements are the children of $R$;
         \item\label{regular d} for each set $R$, for every $\varepsilon>0$ and for every $x\in R$, there exists a descendant $R_{x,\varepsilon}\in \children^\ast(R)$ containing $x$ and with $\diam_S\tonde{R_{x,\varepsilon}}<\varepsilon$;
         \item \label{regular R} the quantity $\displaystyle\mathcal M(S):=\sup_{R \in \admis} \sup_{R' \in \children (R)} \tfrac{\rho(R)}{\rho(R')}$ is bounded;
           \item if $R^*$ is defined as in \eqref{Def-enlargedset}, then \label{regular CZ} $\displaystyle\sup_{R \in \admis} \tfrac{\rho(R^\ast)}{\rho(R)}<\infty$.
        \end{enumerate}
\end{theorem}

Consequently, the abstract Hardy--$\mathrm{BMO}$ and CZ theory developed in \cite{CMV26} becomes available in the present setting and the associated atomic Hardy space \(H^1(S)\) and the corresponding $\BMO(S)$ space can be defined using the admissible sets in \(\admis\).

A natural application concerns the Riesz transforms associated with the distinguished Laplacian on \(S\) defined as \begin{equation*}
    \Delta:=-\sum_{j=0}^{\dimV+\dimZ} X_j^2,
\end{equation*} 
where $\graffe{X_j}_{j=0}^{\dimV+\dimZ}$ is the orthonormal frame of left-invariant vector fields on $S$ introduced above. The distinguished Laplacian takes this form because the left-invariant vector fields $X_j$ are skew-adjoint with respect to the right Haar measure $\rho$, so $\Delta$ is (essentially) self-adjoint on $L^2(S)$.

Let
\[
\Rz_j=\vfX_j\Delta^{-1/2},\,\,\,j=0,\ldots,\dimV+\dimZ,
\]
be the first-order Riesz transforms associated with the distinguished
Laplacian \(\Delta\).  The weak type \((1,1)\) and \(L^p\)-boundedness of \(\Rz_j\) for every \(1<p<\infty\) were recently established in \cite{LM26}. Thanks to our construction of admissible sets on $S$, we can now investigate Hardy-$\mathrm{BMO}$ endpoint bounds for the Riesz transforms $\Rz_j$ in this setting.  

The following theorem summarizes the endpoint mapping properties proved in this paper.
\begin{theorem}\label{Thm-mainendpoint}
The following statements hold for every \(j=0,\dots,\dimV+\dimZ\).
\begin{enumerate}
    \item \label{H1L1bdd}The Riesz transform \(\Rz_j\) extends to a bounded operator from $H^1(S)$ to $L^1(S)$; and \(\Rz_j^*\) extends to a bounded operator from $L^\infty(S)$ to $\BMO(S)$.
    \item \label{H1L1ubdd}The adjoint Riesz transform \(\Rz_j^*\) is not bounded from \(H^1(S)\) to \(L^1(S)\).
    \item \label{BMOubdd}The Riesz transform \(\Rz_j\) is not bounded from \(L^\infty(S)\) to \(\BMO(S)\).
\end{enumerate}
\end{theorem}

The first assertion follows by combining the concrete CZ structure constructed in this paper with the kernel estimates established in \cite{LM26}. The second assertion is proved by constructing suitable families of \(H^1\)-atoms whose images under the Riesz transforms are not integrable, while the third is a consequence of the second one and the \(H^1\)--\(\BMO\) duality. 
Thus, although both \(\Rz_j\) and \(\Rz_j^*\) are bounded on \(L^p(S)\) for every \(1<p<\infty\), their endpoint mapping properties are asymmetric. A similar endpoint asymmetry was observed by Sj\"ogren and Vallarino on the $ax+b$ group $G_2=\RR^2\rtimes\RR^+$: the first-order Riesz transforms are bounded from $H^1(G_2)$ to $L^1(G_2)$, whereas none of their adjoints is \cite[Theorems 3.2, 4.2 and 5.2]{SV08}. Martini subsequently extended the latter negative result to
$G_d=\RR^d\rtimes\RR^+$ for every $d\geq1$, proving that no adjoint first-order Riesz transform is bounded from $H^1(G_d)$ to $L^1(G_d)$; see \cite[Proposition 7.1]{M23}.

\subsection*{Notation}
We denote the set consisting of nonnegative integers by $\NN$, while \(\RR^{+}\) is the set of nonnegative real numbers. The symbol \(\bigsqcup\) denotes the disjoint union. For nonnegative quantities \(A\) and \(B\), we write \(A\lesssim B\) if there exists a constant \(C>0\) such that \(A\le C B\). We write \(A\asymp B\) if both \(A\lesssim B\) and \(B\lesssim A\). If the implicit constant depends on a parameter \(\varepsilon\), we write \(A\lesssim_\varepsilon B\).

\section{Preliminaries} 

\subsection{H-type groups}
\label{subsection-$H$-type}
Let $\lie{n}$ be a 2-step stratified Lie algebra, equipped with an inner product $\langle\argomento,\argomento\rangle$, and denote by $|\argomento|$ the corresponding norm. We assume that $\lie{n}$ decomposes as the orthogonal direct sum
\[
\lie{n}=\lie{v}\oplus\lie{z},
\]
where $\lie{z}=[\lie{n},\lie{n}]$ and $[\lie{n},\lie{z}]=\graffe{0}$. We let $\dimV=\dim\lie{v}$ and $\dimZ=\dim\lie{z}$, and identify the horizontal and vertical layers with the Euclidean spaces $\mathbb{R}^{\dimV}$ and $\mathbb{R}^{\dimZ}$, respectively.
 For any $z \in \lie{z}$, we denote by $J_z$ the skew-symmetric endomorphism of $\lie{v}$ defined by
\begin{equation*}
\langle J_z x,y\rangle=\langle z,[x,y]\rangle\quad\forall x,y\in\lie{v},\quad\forall z\in\lie{z}.
\end{equation*}
We say that $\lie{n}$ is of Heisenberg type ($H$-type for short) if
\begin{equation}\label{eq:Htype}
    |J_{z}x|=|z||x|\quad\forall x\in\lie{v},\quad\forall z\in\lie{z}.
\end{equation}
For later use in Section \ref{Sec-endpointRiesz}, we record the following consequence of the $H$-type condition.
\begin{lemma}
For all \(x,x'\in\mathfrak v\),
\begin{equation}\label{Relation-Htype1}
    \abs{\langle x,x'\rangle}^{2}+\abs{[x,x']}^{2}\leq |x|^{2}|x'|^{2}.
\end{equation}
\end{lemma}
\begin{proof} The assertion is trivial when \(x'=0\). Suppose that \(x'\neq 0\), and define
\[
    x_{\perp} \defeq x-\left\langle x,\frac{x'}{|x'|}\right\rangle\frac{x'}{|x'|}.
\]
Then,
\[
    [x,x']=[x_{\perp},x'],\,\,\,|x_{\perp}|^{2} =|x|^{2}-\frac{\abs{\langle x,x'\rangle}^{2}}{|x'|^{2}}.
\]
Using the definition of \(J_z\) and the $H$-type condition, we obtain
\begin{align*}
    \abs{[x,x']}&=\abs{[x_{\perp},x']} 
    =\sup_{\substack{z\in\mathfrak z\\ |z|=1}}\abs{\left\langle z,[x_{\perp},x']\right\rangle}
    =\sup_{\substack{z\in\mathfrak z\\ |z|=1}} \abs{\left\langle J_zx_{\perp},x'\right\rangle}\\
    &\leq\sup_{\substack{z\in\mathfrak z\\ |z|=1}}|J_zx_{\perp}|\,|x'|
    =|x_{\perp}|\,|x'|.
\end{align*}
Consequently,
\[
    \abs{[x,x']}^{2}\leq|x_{\perp}|^{2}|x'|^{2}= |x|^{2}|x'|^{2}-\abs{\langle x,x'\rangle}^{2},
\]
which proves \eqref{Relation-Htype1}.
\end{proof}

Let $N$ be an $H$-type group. Via the exponential map, we identify the group $N$ with its Lie algebra $\lie{n}$ by 
\[
\begin{array}{rcl}
	\lie{v}\times\lie{z}&\to&N\\
	(x,z)&\mapsto&\exp(x+z).
\end{array}
\]	
The group law in $N$ is given by
\begin{equation}
\label{eq-Ngrouplaw}
(x,z) \Nprod (x',z')= \left( x+x',z+z'+\frac{1}{2}[x,x']\right)
\end{equation}
for all $(x,z), (x',z') \in N$, with $(0,0)$ as the identity element, and moreover
\[
(x,z)^{-1} = (-x,-z).
\]
The group $N$ is unimodular, and the measure $\dd \mu(x,z) :=\dd x  \dd z$ is a left and right Haar measure on $N$.

Let $\RR^+ =(0,+\infty).$ The group $N$ is equipped with the natural family of automorphic dilations $(\dil_t)_{t \in \RR^+}$ defined in \eqref{Def-dilation}.
Let $\dimHomog\defeq\tfrac{\dimV+2\dimZ}2$ denote the homogeneous dimension of $N$ associated with the dilations \((\dil_t)_{t\in \RR^+}\).

On $N$, we introduce two left-invariant metrics, induced by the inner product $\langle\argomento,\argomento\rangle$ on $\lie{n}$. The first is the Carnot--Carath\'eodory metric associated to the horizontal distribution spanned by $\lie{v}$, equipped with the restriction of $\langle\argomento,\argomento\rangle$ to $\lie{v}$; we refer to it as the \emph{sub-Riemannian distance} and denote it by $\subriemdist$. The second is the geodesic distance induced by the Riemannian metric obtained by extending $\langle\argomento,\argomento\rangle$ to the whole group via left translation; we refer to it as the \emph{Riemannian distance} and denote it by $\riemdist$. 

Since both the Carnot--Carath\'eodory distance $\subriemdist$ and the Riemannian distance $\riemdist$ are left-invariant, we have
\begin{equation*}
    \subriemdist(n_1,n_2) = \subriemdist\big((0,0),\, n_1^{-1}n_2\big), \ \ 
    \riemdist(n_1,n_2) = \riemdist\big((0,0),\, n_1^{-1}n_2\big) 
    \ \  \forall\, n_1, n_2 \in N.
\end{equation*}
The following explicit asymptotic estimates hold: for every $n = (x,z) \in N$,
\begin{equation}
\label{Equivalence-KoranyiSR}
    \subriemdist\big((0,0),(x,z)\big) \asymp |x| + |z|^{1/2},
\end{equation}
\begin{equation}
\label{eq-riemannianmetricdefinition}
    \riemdist\big((0,0),(x,z)\big) \asymp 
    \begin{cases}
        |x| + |z|^{1/2}, & \text{if } \max\{|x|,|z|\} \geq 1, \\[2pt]
        |x| + |z|, & \text{otherwise.}
    \end{cases}
\end{equation}
Moreover, $\subriemdist$ is homogeneous with respect to the dilations defined above, in the sense that
\begin{equation}
\label{dilation}
\subriemdist(\dil_t n_1, \dil_t n_2) = t^{1/2}  \subriemdist(n_1,n_2) \quad \forall n_1,n_2 \in N, \quad \forall t\in \RR^+.
\end{equation}
If we define for every $r>0$ and $n_0 \in N$
\[ 
\subriemball {n_0} {r}:=\graffe{n \in N \tc \subriemdist\tonde{n_0, n}<r}, \,\,\, \riemball {n_0} {r}:=\graffe{n \in N \tc \riemdist\tonde{n_0, n}<r},
\]
we then have
\begin{align*}
   \mu\tonde{\subriemball {n_0} {r}} &\asymp r^{\dimV + 2\dimZ} = r^{2\dimHomog}, \quad \forall r >0;\\
    \mu\tonde{\riemball {n_0} {r}} &\asymp \begin{cases} 
        r^{\dimV + \dimZ}, \quad &\mathrm{if \ } r<1; \\
        r^{\dimV + 2\dimZ}, &\mathrm{otherwise.}
    \end{cases}
\end{align*}

We shall use the following ball-box estimates on \(N\). According to \eqref{Equivalence-KoranyiSR}, there exist constants
\(\alpha_1,\beta_1>0\) such that for every \(r>0\),
\begin{equation}\label{Eq-ball-box-N}
\horiball{0}{\alpha_1r}\times \vertball{0}{\alpha_1r^2} \subseteq \subriemball{(0,0)}{r} \subseteq \horiball{0}{\beta_1 r}\times \vertball{0}{\beta_1 r^2}.
\end{equation}
Here \(\horiball{x}{r}\) and \(\vertball{z}{r}\) denote Euclidean balls in the horizontal and vertical layers, respectively. 
\subsection{Damek--Ricci spaces}
\label{subsection-Damek–Ricci spaces}
Damek--Ricci spaces form a family of noncompact, simply connected, solvable Lie groups that generalize the rank-one symmetric spaces of noncompact type. Each Damek--Ricci space is obtained from an $H$-type group \(N\) through a one-dimensional extension.

More precisely, let \(N\) be an $H$-type group with Lie algebra 
\(
\lie{n}=\lie{v}\oplus\lie{z},
\)
and let \(\lie{u}\) be a one-dimensional Lie algebra endowed with an inner product and spanned by a unit vector \(H\). We define the Lie algebra
\[
\lie{s}=\lie{n}\rtimes_\varphi\lie{u}
\]
as the semidirect product of the two Lie algebras, with \(\varphi:\lie{u}\to \mathrm{Der}(\lie{n})\) (the derivations of \(\lie{n}\)) given by
\[
\varphi(H)(X)=\frac{1}{2}X \qquad \forall X\in\lie{v}, \qquad \varphi(H)(Z)=Z \qquad \forall Z\in\lie{z}.
\]

The inner products on \(\lie{n}\) and \(\lie{u}\) are extended to \(\lie{s}\) by imposing that the two subspaces are orthogonal. Let \(N\), \(U\), and \(S\) denote the connected and simply connected Lie groups associated with \(\lie{n}\), \(\lie{u}\), and \(\lie{s}\) respectively. The group \(S\) is called the harmonic extension of the $H$-type group \(N\), and is referred to as a \emph{Damek--Ricci space}.

Using the exponential map, every element \(s\in S\) can be uniquely represented in the form
\[
s =\exp(X+Z)\exp(uH),
\]
with \(X\in\lie{v}\), \(Z\in\lie{z}\), and \(u\in\RR\) ($uH\in \lie u$). Thus, the map
\begin{align*}
\lie{v}\times\lie{z}\times\RR&\longrightarrow S \\ 
(X,Z,u)&\longmapsto \exp(X+Z)\exp(uH)
\end{align*}
provides global coordinates on \(S\). Consequently, \(S\) can be viewed as the semidirect product
\[
S=N\rtimes\mathbb{R},
\]
where \(\mathbb{R}\) acts on \(N\) through the automorphic dilations of the group defined in \eqref{Def-dilation}
\begin{equation*}
\dil_{\ee^u}(x,z)=\bigl(\ee^{u/2}x,\ee^u z\bigr),
\qquad u\in\mathbb{R}.
\end{equation*}

For \(\mathbf{x}=(x,z,u)\) and \(\mathbf{y}=(x',z',u')\), with
\(x,x'\in\mathbb{R}^{\dimV}\), \(z,z'\in\mathbb{R}^{\dimZ}\), and
\(u,u'\in\mathbb{R}\), the group law on \(S\) is given by
\begin{align*}
(x,z,u)\Sprod(x',z',u') &= \bigl(x+\ee^{u/2}x', \ z+\ee^u z'+\tfrac12\ee^{u/2}[x,x'], \ u+u'\bigr)\ \\ 
&= \bigl((x,z)\Nprod\dil_{\ee^u}(x',z'), \ u+u'\bigr).
\end{align*}
The identity element is \(e=(0,0,0)\). We endow $S$ with the left-invariant Riemannian metric obtained by extending the inner product on $\lie{s}$ at the identity element $e$ by left translations.

We denote by $\damekriccidist(\argomento,\argomento)$ the associated geodesic distance on $S$, which satisfies
\[
\damekriccidist(\mathbf{x},\mathbf{y}) = |\mathbf{x}^{-1}\mathbf{y}|,
\]
where $|\mathbf{w}| = \damekriccidist(\mathbf{w},e)$. An explicit formula for $|\mathbf{x}|$ is
\begin{equation}\label{Eq-distanceinS}
\begin{split}
\cosh^2\frac{|\mathbf{x}|}{2}
= \left( \cosh\frac{u}{2} + \frac12 \ee^{-\frac{u}{2}} \Bigl|\frac{x}{2}\Bigr|^2 \right)^2 + \frac14 \ee^{-u} |z|^2.
\end{split}
\end{equation}
The ball of radius $r>0$ centred at $\mathbf{x}$ in $S$ is denoted by
\[
\damekricciball{\mathbf{x}}{r} := \{ \mathbf{y}\in S : \damekriccidist(\mathbf{x},\mathbf{y}) < r \}.
\]
Letting $\rho$ be the right-invariant Haar measure on $S$, given by \begin{equation} \label{eq-damekriccimeasure}
\dd \rho(x,z,u)=\dd x\dd z\dd u,
\end{equation}
the measure of the ball of radius $r$ and centered at $e$ is
$$
\rho(\damekricciball{e}{r})\asymp \begin{cases}
    r^{\dimV+\dimZ+1} \quad &\mathrm{if \ } r\leq 1\\
    \ee^{{\dimHomog}r}\quad &\mathrm{otherwise}.
\end{cases}
$$
Since $\rho$ is right-invariant, it satisfies 
\[
\rho\bigl( A \Sprod g \bigr) = \rho(A)
\quad\text{and}\quad 
\rho\bigl( g \Sprod A \bigr) = \modular(g)^{-1}\rho(A)
\]
for every measurable set $A\subseteq S$ and every $g\in S$, where $\modular$ is the \emph{modular function} of the Damek–Ricci group $S$, given by
\begin{equation}
\label{eq-modular}
\modular(x,z,u) = \ee^{-\dimHomog u}.
\end{equation}

A fact that we will need later is that for every compact set $K \subseteq S$, the Damek--Ricci distance and the Euclidean distance are equivalent when restricted to $K$; see \cite[Lemma 2.2]{CCM19}.

\subsubsection{\rm{\textbf{The $ax+b$ group}}} 
A particular case of a Damek--Ricci space is obtained by taking \(\dimV=0\). In this case, the \(N\) component consists only of the vertical layer \(Z\) and is therefore abelian. The resulting solvable group is an \(ax+b\) group.

Let \(G_{\dimZ}=\RR^{\dimZ}\rtimes\RR\) be an \(ax+b\) group, whose elements we write as \((y,v)\in\RR^{\dimZ}\rtimes\RR\).
Endowing $G_{\dimZ}$ with the left-invariant metric described above for the general case $S$, we denote by \(\axbdist\bigl((y,v),(y',v')\bigr)\) the associated geodesic distance on \(G_{\dimZ}\), and write
\[
|(y,v)|_{ax+b}\defeq\axbdist\bigl((y,v),(0,0)\bigr).
\]
This quantity satisfies the following equivalent equalities
\[
\cosh ^2\frac{|(y,v)|_{ax+b}}{2}=\cosh^2 \frac v2+\frac14\ee^{-v}|y|^2, \quad \cosh {|(y,v)|_{ax+b}}=\cosh v +\frac12\ee^{-v}|y|^2
\]
and also the following coordinate estimates (see for instance \cite[(ii) in Lemma 3.9]{MP25}):
\begin{equation}\label{Eq-axbestimates}
|v|\leq |(y,v)|_{ax+b},\,\,\,|y|\leq\sinh |(y,v)|_{ax+b}.
\end{equation}

Going back to general Damek--Ricci spaces, we obtain the following lemma. 
\begin{lemma}\label{Lem-coordinateestimates}
For every $\mathbf{x}=(x,z,u)\in S$, one has
\[
|x|\leq 2\sinh\frac{|\mathbf{x}|}{2},
\,\,\,
|z|\leq \sinh|\mathbf{x}|,
\,\,\,
|u|\leq |\mathbf{x}|.
\]
\end{lemma}
\begin{proof}
 By \eqref{Eq-distanceinS} and the duplication formula for $\cosh$ ($\cosh(2a) = 2\cosh^2 a -1 \ \forall a \in \RR$):
	$$
	\begin{aligned}
		\cosh |\mathbf{x}|&=2\left(\cosh\frac{u}{2}+\frac{1}{2}\mathrm{e}^{-\frac{u}{2}}\left|\frac{x}{2}\right|^2\right)^2+\frac{1}{2}\mathrm{e}^{-u}|z|^{2}-1\\
		&\geq2\cosh^2\frac{u}{2}-1+\frac{1}{2}\mathrm{e}^{-u}|z|^{2}=\cosh u+\frac{1}{2}\mathrm{e}^{-u}|z|^{2}\\
        &=\cosh \quadre{\axbdist((z,u),(0,0))},
	\end{aligned}
	$$
    by identifying the subset $\graffe{(0,z,u) \ : \ z\in \RR^{\dimZ}, \ u \in \RR}$ with the $ax+b$ group $\RR^{\dimZ}\rtimes \RR$.
    Hence, \(|\mathbf{x}|\geq \axbdist((z,u),(0,0))\).
    Also, by \eqref{Eq-distanceinS},
	$$
	\begin{aligned}
		\cosh \frac{|\mathbf{x}|}{2}&\geq \cosh\frac{u}{2}+\frac{1}{2}\mathrm{e}^{-\frac{u}{2}}\left|\frac{x}{2}\right|^2=\cosh \axbdist\tonde{\tonde{\frac{x}{2},\frac{u}{2}},(0,0)},
	\end{aligned}
	$$
    by identifying the subset $\graffe{(x,0,u) \ : \ x\in \RR^{\dimV}, \ u \in \RR}$ with the $ax+b$ group $\RR^{\dimV}\rtimes \RR$.
    Hence, \(|\mathbf{x}|\geq 2\axbdist((\frac{x}{2},\frac u2),(0,0))\).
	By applying \eqref{Eq-axbestimates}, we get 
	$$
	\begin{aligned}
		|x|&\leq
        2\sinh \axbdist\tonde{\tonde{\frac{x}{2},\frac{u}{2}},(0,0)}\leq2\sinh \frac{|\mathbf{x}|}{2},\\
		|z|&\leq\sinh \axbdist((z,u),(0,0))\leq\sinh |\mathbf{x}|,\\
		|u|&\leq \axbdist((z,u),(0,0))\leq |\mathbf{x}|,
	\end{aligned}
    $$ as required.
\end{proof}

\section{Riemannian and Sub-Riemannian cubes on \texorpdfstring{$H$}{H}-type groups}\label{Sec-dyadicinN}
Our aim in this section is to construct two types of partitions of \(N\): one adapted to large scales and another better suited to small scales, reflecting the fact that, according to \eqref{eq-riemannianmetricdefinition}, the Riemannian metric is equivalent to the sub-Riemannian metric at large scales and to the Euclidean metric at small scales. For later use on Damek--Ricci spaces, our construction must allow the transition scale between the two geometries to be chosen flexibly. The definition mimics the dyadic partition of $\RR^n$. We follow the construction of a tiling of a $H$-type group $N$ due to Strichartz \cite[§2]{S92}; we refer the reader there for the details.

In order to apply Strichartz's construction, one would need the bilinear map appearing in the group law \eqref{eq-Ngrouplaw} at the second layer, namely $\frac{1}{2}\quadre{\argomento,\argomento}$, to have integer coefficients when expressed with respect to suitable bases of $\lie{v}$ and $\lie{z}$. By \cite{CD02}, there exists an orthonormal basis for which the Lie bracket $\quadre{\argomento,\argomento}$ has integer coefficients, and, as we did above, we shall identify $\lie{v}$ and $\lie{z}$ with $\RR^{\dimV}$ and $\RR^{\dimZ}$ respectively, using those bases, endowing the group with the Lebesgue measure. We shall therefore apply Strichartz's construction from \cite[§2]{S92} to the coordinates $(x, 2z)$ instead of $(x,z).$

\subsection{Sub-Riemannian tilings of \texorpdfstring{$N$}{N}}
Let us define, for $j =1,\dots,\dimZ$, the bounded measurable functions on $\RR^{\dimV}$ 
\begin{equation}\label{For_Aj}
A_j(x) =\sum_{n=1}^{\infty} \frac{1}{4^n} \bigg[\lfloor2^n x\rfloor_{\text{mod 2}},\graffe{2^n x}\bigg]_j, \quad x \in \RR^{{\dimV}},
\end{equation}
where $\quadre{\argomento,\argomento}_j$ denote the $j$-th component of the Lie bracket, and $\lfloor \argomento\rfloor$ and $\graffe{\argomento}$ are respectively the integer part and the mantissa functions.
\begin{definition}
   We define the \emph{fundamental tile} at level $0$ as
\begin{equation}\label{Def-T}
T := \graffe{(x,z) \in N \ : \ x \in [0,1)^{\dimV}, \ 2z_j \in \left[A_j(x), A_j(x)+1\right), \ j = 1,\dots,\dimZ}.
\end{equation} 
\end{definition}

The boundedness of the $A_j$'s gives us the following metric property of $T$ .
\begin{proposition}\label{subriemannian balls}
   There exist two positive numbers $0<c_N \leq C_N < \infty$ and a point $n_0^0\in N$ such that 
   \[
    \subriemball{n_0^0}{c_N} \subseteq T \subseteq \subriemball{n_0^0}{C_N} .
   \]
\end{proposition}
\begin{proof}
$T$ has an interior of full $\mu$-measure (\cite[Corollary 2.3.]{S92}), and hence there exist a point $n_0^0\in T$ and a positive number $c_N$ such that $\subriemball{n_0^0}{c_N} \subseteq T$. Also, since $A_j \in L^\infty(\RR^{\dimV})$ for every $j \in \graffe{1,...,\dimZ}$, the set $T$ is bounded and hence there exists a positive constant $C_N >c_N$ such that $T\subseteq \subriemball{n_0^0}{C_N}$.
\end{proof}

\begin{definition}
    A set $A \subseteq N$ is called \emph{self-similar} if there exists a finite set of points $\Lambda\subseteq N$ such that $$\dil_4A=\bigsqcup_{\lambda \in \Lambda}\lambda \Nprod A.$$
\end{definition}
\noindent The definition above is equivalent to the fact that $A=\diltile{-1}{\bigsqcup_{\lambda \in \Lambda}\lambda \Nprod A}.$
\begin{remark}
    The factor $4$ appearing in the dilation $\dil_{4}$ reflects the shape of the dilation on the group $N$. Indeed, since $\dil_t(x,z) = (t^{1/2}x, tz)$, choosing $t=4$ scales the horizontal components by $2$ and the vertical ones by $4$. In particular, the self-similarity condition means that $A$ decomposes into $|\Lambda|$ translated copies of itself at half the horizontal scale.
\end{remark}
\begin{definition}
    Let $T$ be as in \eqref{Def-T}. We define the \emph{dyadic sub-Riemannian cubes} at level $k\in \ZZ$ as
\begin{equation*}
    \Q^k_\gamma:=\gamma\Nprod\dil_{4^{-k}}(T), \quad k\in \ZZ, \gamma\in N.
\end{equation*}
\end{definition}
The term \emph{dyadic} reflects the analogy with the classical dyadic decomposition of \(\RR^n\): each scale is obtained from a fixed one by scaling by powers of $2$, following the geometry of the group $N$. The construction must nonetheless be adapted to the non-Euclidean structure of \(N\), taking into account its group structure and the associated dilations. We will often drop the term \emph{dyadic}, however, and simply refer to these sets as sub-Riemannian cubes.

\begin{theorem}\label{Thm-Christcubes}
Let $c_N $ and $C_N$ be the positive constants appearing in Proposition \ref{subriemannian balls}. The following properties hold:
	\begin{enumerate}
		\item\label{teo3.5-item1} for every $k \in \ZZ$ and $\gamma \in N$, there exists a partition of $N$ of sub-Riemannian cubes of level $k$ where $\mathcal{\Q}^k_\gamma$ is one of the cubes of the partition:
    \begin{equation}
    \label{eq-partitionofN}
        N = \bigsqcup_{\lambda \in\gamma\Nprod \dil_{4^{-k}}(\Gamma)} \mathcal{\Q}^k_\lambda,
    \end{equation}
    for \(\Gamma=\ZZ^{\dimV}\times \tonde{\tfrac 1 2 \ZZ }^{\dimZ}\);
            \item \label{teo3.5-item2} for every $k \in \ZZ$ and every $\gamma \in N$, $$\mu(\Q^k_\gamma)=2^{-k(\dimV+2\dimZ)-\dimZ};$$
		\item \label{teo3.5-item3} $T$ is \emph{self-similar} and each set $\Q_{\gamma}^{k}$ is the disjoint union of $2^{\dimV +2\dimZ}$ sets of the form $\Q_{\lambda}^{k+1}$, which we shall call \emph{sub-Riemannian children} of $\Q^k_\gamma$:
        \begin{equation}
        \label{eq-partitionsubriemanniansubes}
             \Q_\gamma^k =\bigsqcup_{\beta \in \Gamma_{k+1}}\Q^{k+1}_{\gamma \Nprod \beta},
        \end{equation}
        with $\Gamma_k:=\graffe{0,\frac 1 {2^k}}^{\dimV}\times\graffe{0,\frac 1 {2\cdot 4^k},\frac 2 {2\cdot 4^k}, \frac 3 {2\cdot 4^k} }^{\dimZ};$
        \item \label{teo3.5-item4}for every $\Q^k_\gamma$ there exists $n_{\gamma}^{k}$ in $N$ such that 
        \begin{equation}\label{Equ_inclusion}
        \subriemball{n^k_\gamma}{2^{-k} c_N}  \subseteq  \Q^k_\gamma \subseteq  \subriemball{n^k_\gamma}{2^{-k} C_N};
        \end{equation}
        \item\label{teo3.5-item5} we have the following estimates on the sub-Riemannian diameter of $\Q^k_\gamma$: $$c_N 2^{-k+1} \leq \diam_{\mathrm{SR}}(\Q_{\gamma}^{k}) \leq C_N 2^{-k+1}.$$
	\end{enumerate}
\end{theorem}
\begin{proof}
For \eqref{teo3.5-item1}, the case $\Q^k_\gamma$ with $\gamma=(0,0)$ and $k=0$ is already in \cite[Theorem 2.2]{S92}, namely
\begin{equation*}
N = \bigsqcup_{\lambda \in \Gamma}\lambda \Nprod T= \bigsqcup_{\lambda \in \Gamma}\Q^0_\lambda
\end{equation*}
which is a tiling of $N$. As in \cite{S92}, we call it a \emph{stacked tiling}, since tiles with the same base in the $x$-variable have the same shape, and are stacked on top of each other in the $z$-direction. 
If $k\neq 0,$ then 
\begin{align*}
    N=\diltile{-k}{N}= \diltile{-k}{\bigsqcup_{\lambda \in \Gamma}\Q^0_\lambda}=\bigsqcup_{\lambda \in \diltile{-k}\Gamma}\Q^k_\lambda.
\end{align*}
If also $\gamma \neq (0,0)$, we then have
 \[
 N=\gamma\Nprod N=\gamma\Nprod \bigsqcup_{\lambda \in \diltile{-k}\Gamma}\Q^k_\lambda = \bigsqcup_{\lambda \in \gamma \Nprod\diltile{-k}\Gamma}\Q^k_\lambda.
 \] 
Clearly, $\Q^k_\gamma$ is one of the cubes of the partition, since the point $(0,0) \in \diltile{-k}{\Gamma}$ for every $k\in \ZZ.$ We conclude \eqref{teo3.5-item1}.

Fubini's Theorem and the fact that the measure on the Heisenberg-type group $N$ is the Lebesgue measure $\mu$ yields 
$$\mu(T)= 2^{-\dimZ},$$ since each $x$ coordinate ranges over an interval of length $1$ and each $z$ coordinate ranges in an interval of length $\tfrac 1 2.$
Lebesgue measure on $N$ is also left- and right-invariant, so, for every $\gamma \in N,$
$$
\mu(\Q^0_\gamma) = \mu(\gamma \Nprod T) =\mu(T)= 2^{-\dimZ}.
$$
The dilation operator $\dil_{4^{-k}}$ halves the length of the $x$-components $k$ times, and halves the length of the $z$-components $2k$ times. Hence 
$$
\mu(\Q^k_\gamma) = \mu(\gamma \Nprod \diltile{-k}T) =\mu(\diltile{-k}T)= 2^{-\dimV k} 2^{-2\dimZ k}2^{-\dimZ} =2^{-k\tonde{\dimV +2\dimZ}-\dimZ},
$$
which gives \eqref{teo3.5-item2}.

Now we prove \eqref{teo3.5-item3}. First, $T$ is self-similar: thanks to \cite[Lemma 2.1]{S92}, the functions $A_j$ satisfy 
\begin{equation} \label{Aj} 
A_j(x) = \frac{1}{4} A_j(\langle 2x \rangle) + \frac{1}{4}     \left[\lfloor 2x \rfloor, \langle 2x \rangle \right]_j, \quad x \in [0,\tfrac12)^{\dimV},
\end{equation} which implies \cite[Theorem 2.2]{S92}:
    \begin{equation}
    \label{self similarity}
      \diltile{}{T} = \bigsqcup_{\lambda \in \Gamma_0}\lambda\Nprod T =   \bigsqcup_{\lambda \in \Gamma_0}\Q^0_\lambda
    \end{equation}
    or equivalently
\begin{equation}
\label{self similarity 2}
T =\bigsqcup_{\lambda \in \diltile{-1}{\Gamma_0}}\lambda\Nprod \dil_{4^{-1}}\tonde{ T},
\end{equation}
where
$$\Gamma_0:=\{0,1\}^{\dimV}\times \left\{0,\frac{1}{2},1, 
\frac{3}{2}\right\}^{\dimZ}.$$ In particular, the fact that the union above is disjoint comes from the inclusion $\Gamma_0\subseteq \Gamma$.

If we now take an arbitrary $\Q^k_\gamma$, denoting $\Gamma_k:=\dil_{4^{-k}}\Gamma_0$ ($\subseteq \diltile{-k}{\Gamma}$), we have that
$$
\begin{aligned}
    \Q_\gamma^k &=\gamma\Nprod \dil_{4^{-k}}( T) =\gamma \Nprod \dil_{4^{-k}}\tonde {\bigsqcup_{\lambda\in \Gamma_0}\dil_{4^{-1}}\tonde{\lambda \Nprod T}} \\
    &= \bigsqcup_{\lambda \in \Gamma_{0}}\gamma \Nprod \dil_{4^{-k-1}} \lambda\Nprod  \dil_{4^{-k-1}}(T)=\bigsqcup_{\beta \in \Gamma_{k+1}}\Q^{k+1}_{\gamma \Nprod \beta}.
\end{aligned}$$
This subdivision process produces $\abs{\Gamma_{k+1}}=\abs{\Gamma_0}= 2^{\dimV+2\dimZ}$ \emph{sub-Riemannian children} of $\Q_\gamma^k$. Thus, \eqref{teo3.5-item3} is proved.

Next we prove \eqref{Equ_inclusion}. From Proposition \ref{subriemannian balls}, we have
	$$
	\subriemball{n_0^0}{c_N} \subseteq T\subseteq \subriemball{n_0^0}{C_N}.
	$$
    Moreover, since the sub-Riemannian distance on $N$ is left-invariant and homogeneous by \eqref{dilation}, we obtain that $\forall r>0, \ \gamma_1,\gamma_2 \in N, \ k\in \ZZ$
	$$
	\dil_{4^{-k}}\subriemball{n^0_0}{r}=\subriemball{\diltile{-k}{n_0^0}}{2^{-k}r}, \quad \gamma_1\Nprod\subriemball{\gamma_2}{r}=\subriemball{\gamma_1 \Nprod\gamma_2}{r},
   $$ 
    which means that the metric properties of $T$ transfer to every cube $\Q^k_\gamma = \gamma \Nprod \dil_{4^{-k}}(T)$. Indeed, denoting by $n^k_\gamma := \gamma\Nprod\dil_{4^{-k}}( n_0^0)$, the inclusions above yield \eqref{Equ_inclusion}.

	\eqref{teo3.5-item5} can be obtained easily from \eqref{Equ_inclusion}.
\end{proof}
   
We now introduce a definition that will be used later, which establishes a ``hierarchy'' among the cubes.
\begin{definition}\label{Def-Subriemannian generation}
Let $\Q^k_\gamma$ be a cube of level $k$. If 
$$
\Q^k_\gamma =\bigsqcup_{i=1}^{2^{\dimV + 2\dimZ}}\Q^{k+1}_{\gamma_i}
$$
as in \eqref{teo3.5-item3} of Theorem \ref{Thm-Christcubes}, we denote by $\mathcal{C}^1_{\mathrm{SR}}(\Q^k_\gamma)$ the collection $\graffe{\Q^{k+1}_{\gamma_i}}_{i=1}^{2^{\dimV+2\dimZ}}$ and we call it the \textit{first sub-Riemannian generation} (children) of $\Q^k_\gamma$.

More generally, for $h\geq 2$, we define the \textit{$h$-th sub-Riemannian generation} of $\Q^k_\gamma$ inductively as $$ \mathcal{C}^{h}_{\mathrm{SR}}(\Q^k_\gamma)= \mathcal{C}^1_{\mathrm{SR}}(\mathcal{C}^{h-1}_{\mathrm{SR}}(\Q^k_\gamma)).$$
\end{definition}

\begin{remark}
\label{rmk-christ}
As a sanity check, the sub-Riemannian cubes can in fact be recovered as instances of the abstract dyadic structure of Hyt\"onen and Kairema in \cite{HK12}, which in turn generalizes Christ's construction in \cite{C90}.

If, in particular, we choose, for every $k \in M_0\ZZ$ (for an integer $M_0$ to be specified later), the partition containing the cube $\Q^k_0$ (this is just one possible choice; any other choice works equally well, as long as at every generation the cubes partition those of the previous generation, without any overlapping), we obtain a collection of cubes satisfying the hypotheses of \cite[Theorem 2.2]{HK12} by Hyt\"onen and Kairema, with respect to the sub-Riemannian metric.

Indeed, consider the points $n^k_\gamma$. Since by \eqref{Equ_inclusion} 
\[
\subriemball{n^k_\gamma}{2^{-k}c_N} \subseteq \Q^k_\gamma
\]
and the cubes $\Q^k_\gamma$ are pairwise disjoint, it follows that
\[
\subriemdist(n^k_\gamma, n^k_\lambda) \geq 2\, c_N\, 2^{-k}, \qquad \gamma \neq \lambda.
\]
At the same time, since the cubes $\{\Q^k_\gamma\}_{\gamma\in\diltile{-k}{\Gamma}}$ partition $N$ and by \eqref{Equ_inclusion}
\[
\Q^k_\gamma \subseteq \subriemball{n^k_\gamma}{2^{-k}C_N},
\]
every $x \in N$ satisfies
\[
\min_{\gamma\in \diltile{-k}{\Gamma}} \subriemdist(x, n^k_\gamma) < C_N\, 2^{-k}.
\]
Thus the points $\{n^k_\gamma\}_{k\in \ZZ,\ \gamma \in\diltile{-k}{\Gamma}}$ satisfy condition (2.4) of \cite[Theorem. 2.2]{HK12} with $c_0 = 2c_N$, $C_0 = C_N$, and $\delta = 2^{-M_0}$, with $M_0$ such that hypothesis (2.3) of \cite[Theorem 2.2]{HK12} is satisfied, which is:
\begin{equation}
\label{eq-12delta}
12C_N 2^{-M_0}\leq c_N.
\end{equation}
The only effect of $M_0$ being not $1$ yields to retaining only every \(M_0\)-th generation of cubes. We therefore define \[
\P_\gamma^k\defeq\Q_\gamma^{M_0k}, \qquad z_\gamma^k\defeq n_\gamma^{M_0k}.
\]

The properties verified above do not yet imply, however, that the cubes generated by Hyt\"onen and Kairema's theorem coincide with our cubes \(\P_\gamma^k\). It remains to verify that the parent relation and the assignment of boundary points in the construction of \cite{HK12} can be chosen so that the resulting half-open cubes are precisely \(\P_\gamma^k\). Since the reference points uniquely determine the interiors and closures of the cubes, and the boundary of the half-open cubes is chosen arbitrarily so as to form a partition, the latter point is not a concern, and it therefore remains to check the former: since a cube is determined by the points of its children, it suffices to verify that the partial order induced by inclusion among our cubes,
\[
(k+1,\beta)\le(k,\gamma) \Longleftrightarrow \P_\beta^{k+1}\subseteq\P_\gamma^k,
\]
coincides with the partial order constructed in \cite[Lemma 2.10]{HK12}, or, in other words, that this order relation satisfies the assumptions of that lemma. We now verify these assumptions one by one.
\begin{itemize}
\item If
\[
\subriemdist(z_{\beta}^{k+1},z_{\gamma}^{k})<c_N\delta^k,
\]
then \(z_{\beta}^{k+1}\in\P_\gamma^k\), by the inner-ball inclusion in \eqref{Equ_inclusion}. Since the partitions are nested, the unique cube of generation \(k+1\) containing \(z_{\beta}^{k+1}\), namely
\(\P_\beta^{k+1}\), is contained in \(\P_\gamma^k\).
Thus \((k+1,\beta)\le(k,\gamma)\).
\item If \((k+1,\beta)\le(k,\gamma)\), then \(z_\beta^{k+1}\in\P_\gamma^k\), and hence the outer-ball inclusion in \eqref{Equ_inclusion} gives
\[
\subriemdist(z_\beta^{k+1},z_\gamma^k)
<C_N\delta^k.
\]
\item Every cube of generation \(k+1\) is contained in exactly one cube of
generation \(k\), since the partitions are nested.
\item Every cube in the family has exactly \(2^{M_0(\dimV+2\dimZ)}\) children, so in particular the number of children is uniformly bounded.
\item For \(\ell\ge k\), \((\ell,\beta)\le(k,\gamma)\) if and only if the inclusion
\(\P_\beta^\ell\subseteq\P_\gamma^k\) can be realised through the unique chain of intermediate ancestors.
\end{itemize}
Hence, since the closure of each cube depends only on its children (and therefore only on elements that are smaller with respect to the order relation), the cubes $\P_\gamma^k$ are precisely those generated by the construction of \cite[Theorem 2.2]{HK12}.

Therefore, in addition to the properties listed in Theorem \ref{Thm-Christcubes}, the following also hold (for the fixed family of partitions of nested cubes):
\begin{enumerate}
    \item if $l \geq k$, then either $\P^l_{\gamma'} \subseteq \P^k_{\gamma}$ or $\P^l_{\gamma'} \cap \P^k_{\gamma} = \emptyset$;
    \item for each $(k,\gamma)$ and each $l < k$ there is a unique $\gamma'$ such that $\P^k_{\gamma} \subseteq \P^l_{\gamma'}$;
    \item there exist $C,\eta>0$ such that \cite[(3.6) of Theorem 11]{C90} holds for all $t>0$
    \begin{equation}
    \label{eq-small boundaries}
    \mu\Big(\Big\{n \in \Q^k_\gamma \ : \ \dist_{\mathrm{SR}}(n, N\setminus \Q^k_\gamma) \leq t\, 2^{-k}\Big\}\Big) \leq C\, t^\eta\mu(\Q_\gamma^k).
    \end{equation}
\end{enumerate}
The first two properties are precisely those established in \cite[Theorem 2.2]{HK12}.
The third property follows from \cite[Lemma 5.12]{HK12}. Indeed, since for every $\Q^k_\gamma$ there exists a nested partition satisfying $\P^h_\gamma = \Q^k_\gamma$, up to rescaling the distance $\subriemdist$ by a factor $2^{-k}$, and since the measure we consider on $N$ is doubling, \cite[Lemma 5.12]{HK12} applies and yields Lemma 17 of Christ \cite{C90}. This in turn gives \cite[(3.6) of Theorem 11]{C90}, which is precisely \eqref{eq-small boundaries}.
\end{remark}

\subsection{Riemannian tilings of \texorpdfstring{$N$}{N}}

Starting from a given cube $\Q^k_\gamma$ of level $k$, we now describe how to construct its successive generations of \emph{Riemannian children}.

Since every cube $\Q^k_\gamma$ is the left translated and dilated of $T$, it suffices to describe the construction on the fundamental tile $T$, dilate it via $\dil_{4^{-k}}$ and then translate it by $\gamma$.
Recall that the $2^{\dimV+2\dimZ}$ sub-Riemannian children of $T$ are obtained by subdividing by a factor of $2$ in the $x$-components and by a factor of $4$ in the $z$-components. For the Riemannian subdivision, we instead subdivide by a factor of $2$ in both components, producing $2^{\dimV+\dimZ}$ children, while preserving the property that $T$ is partitioned into finitely many tiles at each generation, each with the same measure.

Motivated by the fact that
\begin{equation}
\label{eq-norminfinity}
A_j(x)=\frac{1}{4}A_j(2x), \quad \forall x\in \left[0,\frac{1}{2}\right)^{\dimV}, \quad j=1,\cdots,\dimZ,
\end{equation}
which is a consequence of the definition of $A_j$ in \eqref{For_Aj}, and is stated in \cite[Lemma 2.1]{S92}, we provide the following.
\begin{definition} The \emph{first fundamental Riemannian tile} is
\begin{equation*}
T^{\mathrm R,1} \defeq \left\{ (x,z) \in N :  2x \in [0,1)^{\dimV}, 4 (2z_j) \in \left[A_j(2x), A_j(2x) + 2\right), j = 1,\dots,\dimZ \right\}.
\end{equation*}
More generally, for \(h\in \NN\), the \emph{$h$-th fundamental Riemannian tile} is \label{Def-TRh}
\begin{equation*}
T^{\mathrm R,h}\defeq\left\{(x,z)\in N:2^h x\in[0,1)^{\dimV},4^h(2z_j)\in\left[A_j(2^h x),A_j(2^h x)+2^h\right),j=1,\dots,\dimZ\right\},
\end{equation*}
with the notation $T^{R,0}\defeq T$.
\end{definition}
\begin{remark}
Notice that $\Q_0^h$, namely the sub-Riemannian cube of generation $h$ obtained by dilating the fundamental tile $T$ and having its ``bottom-left'' corner at the origin, is given by
\begin{align}
\label{eq-riemannian dilated tile}
\Q_0^h = \diltile{-h}{T}  = &\big\{(x,z)\in N : 2^h x \in[0,1)^{\dimV}, \\
\notag &\quad \quad \ \ \ 4^h(2z_j)\in\big[A_j(2^h x),A_j(2^h x)+1\big),\ j=1,\dots,\dimZ\big\}.
\end{align}
This expression is similar to that of the fundamental Riemannian tile of generation $h$. The difference lies in the size of the intervals containing the $z$-coordinates. Indeed, in the sub-Riemannian case these intervals have length $2^{-2h-1}$, whereas in the Riemannian case they have length $2^{-h-1}$. Thus, when passing from generation $h$ to generation $h+1$, the length of the intervals in the $z$-variables is divided by $4$ in the sub-Riemannian setting, while it is divided by $2$ in the Riemannian one.
\end{remark}
\begin{remark}
    From \eqref{eq-norminfinity}, we get that 
\[
\norma{A_j}_{L^{\infty}[0,1/2)} = \frac 1 4 \norma{A_j}_{L^{\infty}[0,1)},
\]
and proceeding recursively  
\[
\norma{A_j}_{L^{\infty}[0,2^{-h})} = \frac 1 {4^h} \norma{A_j}_{L^{\infty}[0,1)}.
\]
This yields the fact that the Riemannian tile $T^{\mathrm R,h}$ becomes increasingly with $h$ similar to a product of intervals.
\end{remark}

\begin{definition}
    We define the \emph{dyadic Riemannian cubes} at level $(k,h)\in \ZZ\times \NN$ as:
    \begin{equation}\label{Def-Qkhgamma}
\Q^{k,h}_\gamma\defeq\gamma\Nprod\dil_{4^{-k}}(T^{R,h}), \quad k\in \ZZ, h\in \NN,  \gamma\in N,
\end{equation} 
with the convention that $\Q^{k,0}_\gamma\defeq\Q^k_\gamma.$
\end{definition}
As in the sub-Riemannian case, we will often omit the adjective \emph{dyadic} when referring to Riemannian cubes.

We now show that each Riemannian cube at level \((k,h)\) can be seen as a union of sub-Riemannian cubes at level \(k+h\), with the number of constituent cubes increasing with \(h\). This construction preserves some of the geometric features of the sub-Riemannian cube system, while producing tiles for which the ratio between the length of the interval over which each $x_j$-component ranges (one at a time, with all others fixed) and the length of the interval over which each $z_j$-component ranges (one at a time, with all others fixed) is independent of \(h\). The following proposition explains the particular form of the Riemannian cubes introduced above, whose geometry is chosen precisely to make the decomposition in the following proposition possible.

\begin{proposition} \label{prop-riemannianchristcubes}
 For every $h\in \NN$, $k \in \ZZ$ and $\gamma \in N$, the Riemannian cube $\Q^{k,h}_\gamma$ is the disjoint union of $2^{h\dimZ}$ sub-Riemannian cubes of level $k+h$ stacked in the $z$-direction, namely 
 \begin{equation}\label{eq-q+hgamma}
     \Q^{k,h}_\gamma = \bigsqcup_{\lambda \in \gamma\Nprod \diltile{-k}{\Lambda_h} } \Q^{k+h}_\lambda,
 \end{equation} where \begin{equation}\label{Def-Lambdah}
\Lambda_h\defeq\{0\}^{\dimV}\times\left\{\frac{l}{2\cdot4^h}:l=0,\ldots,2^h-1\right\}^{\dimZ}.
\end{equation}
\end{proposition}
\begin{proof}
We first prove the proposition for the Riemannian fundamental tile $T^{\mathrm R,h}$, and then extend the result to all Riemannian cubes by dilation and translation.

Consider the set $\Lambda_h$ as in \eqref{Def-Lambdah} and consider the sub-Riemannian cube \(\Q_0^h\) in \eqref{eq-riemannian dilated tile}. For $\lambda=(\mathbf 0_{\dimV}, \lambda_1,\dots,\lambda_{\dimZ})\in\Lambda_h,$
the left translation of \(\Q_0^h\) by \(\lambda\) gives
\[
\begin{aligned}
\lambda\Nprod\Q_0^h = \Q^h_\lambda
= \bigg\{  (x,z)\in N: \,\,&2^hx\in[0,1)^{\dimV}, 4^h(2z_j)\in [B_j,B_j+1), \\ 
 & B_j\defeq  2\cdot4^h\cdot \lambda_j+A_j(2^hx)+1,  j=1,\dots,\dimZ \bigg\}
\end{aligned}
\]
and this set is contained in $ T^{\mathrm{R},h}$ from the definition in \eqref{Def-TRh} since $2\cdot4^h\cdot \lambda_j\geq 0 $ and $2\cdot4^h\cdot \lambda_j + 1 \leq 2^h$ for each $j=1,\dots,\dimZ.$ 
Also, \begin{align*}
    \bigcup_{\lambda\in \Lambda_{h}}\lambda\Nprod\Q_0^h= \big\{&(x,z)\in N: 2^hx\in[0,1)^{\dimV},\\ 
&4^h(2z_j)\in \left[ A_j(2^hx), A_j(2^hx)+2^h\right), j=1,\dots,\dimZ \big\},
\end{align*}
which is exactly \(T^{\mathrm R, h}\).

Furthermore, since
\[
\Lambda_h\subseteq\diltile{-h}{\Gamma},
\]
the union of $\Q^h_\lambda$ with $\lambda \in \Lambda_h$ is a disjoint union. 
Therefore, we proved
\begin{equation}
\label{eq-lambdah}
   T^{\mathrm R, h} = \bigsqcup_{\lambda\in\Lambda_h} \Q^h_\lambda, 
\end{equation}
where $\abs{\Lambda_h}= 2^{h\dimZ}.$

 The general case of \eqref{eq-q+hgamma} is obtained straightforwardly by dilating and translating both sides of the equality for the fundamental Riemannian tiles \eqref{eq-lambdah}.
\end{proof}
The proposition above leads to the following theorem, which is analogous to Theorem \ref{Thm-Christcubes}, but for Riemannian cubes.
\begin{theorem}\label{Thm-RiemannianChristcubes}
For every cube $\Q^{k,h}_\gamma$ with $k\in \ZZ$, $h \in \NN$ and $\gamma \in N$, the following properties hold:
	\begin{enumerate}
		\item\label{teo3.11-item1} there exists a partition of $N$ of Riemannian cubes of level $(k,h)$ where ${\Q}^{k,h}_\gamma$ is one of the cubes of the partition:
    \begin{equation}
    \label{eq-riemannianpartitionofN}
        N = \bigsqcup_{\lambda \in\gamma\Nprod \Gamma_{k,h}} {\Q}^{k,h}_\lambda,
    \end{equation}
    where
    \begin{equation}
    \label{eq-gammakhdefinizione}
        \Gamma_{k,h} \defeq \diltile{-k}{\tonde{\tfrac{1}{2^h}\ZZ}^{\dimV} \times \left(\tfrac{1}{2^{h+1}}\ZZ\right)^{\dimZ}};
    \end{equation}
     \item \label{teo3.11-item2}
     $$\mu(\Q^{k,h}_\gamma)=2^{-h(\dimV+\dimZ) - k(\dimV+2\dimZ)-\dimZ};$$
		\item \label{teo3.11-item3}each set $\Q_{\gamma}^{k,h}$ is the disjoint union of $2^{\dimV +\dimZ}$ sets of the form $\Q_{\lambda}^{k,h+1}$, which we shall call \emph{Riemannian children} of $\Q^{k,h}_\gamma$:
        \begin{equation}
            \label{eq-deltapartizione}
            \Q_{\gamma}^{k,h}=\bigsqcup_{\lambda\in \Delta_{k,h}} \Q^{k,h+1}_{\gamma \Nprod \lambda},
        \end{equation}
        with $\Delta_{k,h}:=\diltile{-k}{\graffe{0, \frac{1}{2^{h+1}}}^{\dimV} \times\graffe{0, \frac{1}{2^{h+2}}}^{\dimZ}}.$ 
	\end{enumerate}
\end{theorem}
\begin{proof}
    \eqref{teo3.11-item1} This is a consequence of Proposition \ref{prop-riemannianchristcubes} and (1) of Theorem \ref{Thm-Christcubes}. We prove it just for $\Q^{0,h}_0=T^{\mathrm R,h}$: the conclusion follows after dilation and left translation, as we did previously.

By \eqref{eq-lambdah},
$$ T^{\mathrm R,h} = \bigsqcup_{\lambda\in\Lambda_h}\Q^h_\lambda,$$
and by \eqref{eq-partitionofN}, the cubes $\Q^h_\gamma$ form a tiling of $N$:
$$N = \bigsqcup_{\gamma \in \dil_{4^{-h}}(\Gamma)} \Q^h_\gamma.$$
Since $N$ is partitioned into sub-Riemannian cubes, and a Riemannian cube is a union of sub-Riemannian cubes, it only remains to group the sub-Riemannian cubes of the partition so as to obtain a partition into Riemannian cubes. We notice that if we consider $\Lambda_h$ as in \eqref{Def-Lambdah} and take the product of its elements with the elements of
\begin{equation*}
    \label{eq-lambda0h}
    \Gamma_{0,h}=\tonde{\frac{1}{2^h}\ZZ}^{\dimV} \times\tonde{\frac{1}{2^{h+1}}\ZZ}^{\dimZ},
\end{equation*}
we obtain, for a generic $$\lambda=\big(\mathbf 0_{\dimV},\ \tfrac{\lambda_1}{2\cdot4^h},\dots,\tfrac{\lambda_{\dimZ}}{2\cdot4^h}\big)\in\Lambda_h\quad \mathrm{and}\quad \gamma=\tfrac{1}{2^h}\big(n_1,\dots,n_{\dimV},\ \tfrac{m_1}{2},\dots,\tfrac{m_{\dimZ}}{2}\big)\in\Gamma_{0,h},$$ with $n_i,m_j\in\ZZ$,
$$
\lambda\Nprod\gamma
=\tonde{\tfrac{n_1}{2^h},\dots,\tfrac{n_{\dimV}}{2^h},\ \tfrac{2^hm_1+\lambda_1}{2\cdot4^h},\dots,\tfrac{2^hm_{\dimZ}+\lambda_{\dimZ}}{2\cdot4^h}},
$$
which is the generic element of $\diltile{-h}{\Gamma}$, written in a unique way, since each $\lambda_j$ ranges from $0$ to \(2^h-1\). We have thus shown that
$$
N=\bigsqcup_{\gamma\in\Gamma_{0,h}}\Q^{0,h}_\gamma.
$$

By dilation and left translation, \eqref{eq-riemannianpartitionofN} follows easily.\\
        \eqref{teo3.11-item2} This is a consequence of Proposition \ref{prop-riemannianchristcubes} and (2) of Theorem \ref{Thm-Christcubes}: 
        $$\mu(\Q^{k,h}_\gamma)= 2^{h\dimZ} \mu(\diltile{-k-h}{T})=2^{-h(\dimV+\dimZ) - k(\dimV+2\dimZ)-\dimZ}.$$
        \eqref{teo3.11-item3} It follows from Proposition \ref{prop-riemannianchristcubes} and (3) of Theorem \ref{Thm-Christcubes}. Again, we can prove the statement for $\Q^{0,h}_0=T^{\mathrm{R},h}$, as the conclusion follows from dilation and left translation. 

       Since, by \eqref{eq-lambdah},
$$ T^{\mathrm R,h} = \bigsqcup_{\lambda\in\Lambda_h}\Q^h_\lambda \quad \text{and} \quad T^{\mathrm R,h+1} = \bigsqcup_{\lambda\in\Lambda_{h+1}}\Q^{h+1}_\lambda,$$
and, from \eqref{eq-partitionsubriemanniansubes}, every sub-Riemannian cube of generation $h$ splits as
\begin{equation}
\label{eq-gammabeta}
             \Q_\gamma^h =\bigsqcup_{\beta \in \Gamma_{h+1}}\Q^{h+1}_{\gamma\Nprod \beta},
        \end{equation}
with $\Gamma_h:=\graffe{0,\tfrac{1}{2^h}}^{\dimV}\times\graffe{0,\tfrac{1}{2\cdot4^h},\tfrac{2}{2\cdot4^h},\tfrac{3}{2\cdot4^h}}^{\dimZ}$, combining the two decompositions we obtain a partition of $T^{\mathrm R,h}$ into sub-Riemannian cubes of generation $h+1$:
$$ T^{\mathrm R,h} = \bigsqcup_{\lambda\in\Lambda_h}\bigsqcup_{\beta\in\Gamma_{h+1}}\Q^{h+1}_{\lambda\Nprod\beta}.$$

We now regroup these sub-Riemannian cubes of level $h+1$ into translates of $T^{\mathrm R,h+1}$, as follows. If $T^{\mathrm R,h}$ decomposes into translates of $T^{\mathrm R,h+1}$, then
\begin{equation}
\label{eq-delta0h}
T^{\mathrm R,h } = \bigsqcup_{\delta\in\Delta_{0,h}} \delta\Nprod T^{\mathrm R,h+1 } = \bigsqcup_{\delta\in\Delta_{0,h}}  \bigsqcup_{\gamma\in\Lambda_{h+1}}\Q^{h+1}_{\delta\Nprod\gamma},
\end{equation}
and, comparing this with the decomposition above, we must have that, for every $\lambda\in\Lambda_h$ and $\beta\in\Gamma_{h+1}$, there exist unique $\delta\in\Delta_{0,h}$ and $\gamma\in\Lambda_{h+1}$ such that $\lambda\Nprod\beta=\delta\Nprod\gamma$.

We consider
\begin{equation}
\label{eq-delta0hdef}
\Delta_{0,h}:=\graffe{0,\tfrac{1}{2^{h+1}}}^{\dimV}\times\graffe{0,\tfrac{1}{2^{h+2}}}^{\dimZ},
\end{equation}
whose cardinality is $2^{\dimV+\dimZ}$. Before finding $\delta$ and $\gamma$ explicitly, we check that the cardinalities of $\Lambda_h\times\Gamma_{h+1}$ and $\Delta_{0,h}\times\Lambda_{h+1}$ coincide:
\begin{align*}
&\abs{\Lambda_h}\cdot\abs{\Gamma_{h+1}}=2^{h\dimZ}\cdot2^{\dimV+2\dimZ}=2^{\dimV+(h+2)\dimZ},\\
&\abs{\Delta_{0,h}}\cdot\abs{\Lambda_{h+1}}=2^{\dimV+\dimZ}\cdot2^{(h+1)\dimZ}=2^{\dimV+(h+2)\dimZ}.
\end{align*}
It therefore suffices to show the existence, for every $\lambda \in \Lambda_h$ and $\beta\in \Gamma_{h+1}$, of $\delta\in\Delta_{0,h}$ and $\gamma\in\Lambda_{h+1}$ such that $\lambda\Nprod\beta=\delta\Nprod\gamma$.

Since $\lambda$ and $\gamma$ have vanishing horizontal components, $\Nprod$ reduces to the sum $\RR^{\dimV+\dimZ}$, so
\begin{align*}
\lambda\Nprod\beta=(\beta_1,\dots, \beta_{\dimV},\ \lambda_{\dimV+1}+\beta_{\dimV+1}, \dots, \lambda_{\dimV+\dimZ}+\beta_{\dimV+\dimZ}),\\
\delta\Nprod\gamma=(\delta_1,\dots, \delta_{\dimV},\ \gamma_{\dimV+1}+\delta_{\dimV+1}, \dots, \gamma_{\dimV+\dimZ}+\delta_{\dimV+\dimZ}),
\end{align*}
which forces $\delta_j=\beta_j\in\graffe{0,\tfrac{1}{2^{h+1}}}^{\dimV}$, for every $j\in \graffe{1,\dots,\dimV}.$

As for the vertical components, for $j\in\graffe{\dimV+1,\dots,\dimV+\dimZ}$, writing $\lambda_j=\tfrac{k}{2\cdot4^h}$, $k\in\graffe{0,\dots,2^h-1}$, and $\beta_j=\tfrac{\ell}{2\cdot4^{h+1}}$, $\ell\in\graffe{0,1,2,3}$; then
$$
\lambda_j+\beta_j=\frac{4k+\ell}{2\cdot4^{h+1}}=:\frac{m}{2\cdot4^{h+1}},\qquad m\in\graffe{0,\dots,2^{h+2}-1}.
$$
$m$ can be written uniquely as $m=2^{h+1}\varepsilon+r$ with $\varepsilon\in\graffe{0,1}$, $r\in\graffe{0,\dots,2^{h+1}-1}$, and hence
$$
\delta_j=\frac{\varepsilon}{2^{h+2}}\in\graffe{0,\tfrac{1}{2^{h+2}}},\qquad \gamma_j=\frac{r}{2\cdot4^{h+1}}\in\graffe{0,\dots,\tfrac{2^{h+1}-1}{2\cdot4^{h+1}}},
$$
yielding $\lambda_j+\beta_j=\delta_j+\gamma_j$, and hence the equality we looked for.  

We deduce that $T^{\mathrm R,h}$ decomposes into translates of $T^{\mathrm R,h+1}$ indexed by $\Delta_{0,h}$, as claimed in \eqref{eq-delta0h}.

The general statement \eqref{eq-deltapartizione} follows by dilation and left translation.
\end{proof}
Analogously as above, we define the Riemannian generations of cubes.
\begin{definition}\label{Def-Riemannian generation}
Let $\Q^{k,h}_\gamma$ be a Riemannian cube of level $(k,h)$, where $k\in\ZZ$, $h\in\NN$, and $\gamma\in N$. If 
$$
\Q^{k,h}_\gamma = \bigsqcup_{i=1}^{2^{\dimV+\dimZ}} \Q^{k,h+1}_{\gamma_i}
$$
as in \eqref{eq-deltapartizione}, then we denote by $\mathcal{C}^1_{\mathrm{R}}(\Q^{k,h}_\gamma)$ the collection $\graffe{\Q_{\gamma_i}^{k,h+1}}_{i=1}^{2^{\dimV+\dimZ}}$ and we call it the \textit{first Riemannian generation} (children) of $\Q^{k,h}_\gamma$.

More generally, for $\ell\geq 2$, we define the \emph{$\ell$th Riemannian generation} of $\Q^{k,h}_\gamma$ recursively as $$\mathcal C^\ell_{\mathrm R}(\Q^{k,h}_\gamma)= \mathcal C^1_{\mathrm R}(\mathcal C^{\ell-1}_{\mathrm R}(\Q^{k,h}_\gamma)).$$
\end{definition}
\begin{remark}
     Observe that, when $h=0$, $\Q^{k,0}_\gamma=\Q^k_\gamma$ is a sub-Riemannian cube, and hence we defined above also the Riemannian generations of a sub-Riemannian cube.
\end{remark}
We now show that, when given a cube (either Riemannian or sub-Riemannian), the parameters determining it are unique.
\begin{proposition}
    If for some $k,k' \in \ZZ$, $h,h' \in \NN$ and $\gamma,\gamma' \in N$ we have $$\Q^{k,h}_\gamma = \Q^{k',h'}_{\gamma'},$$ then $$k =k', \quad  h =h', \quad  \gamma =\gamma'.$$
\end{proposition}
\begin{proof}
    We show that the parameters $k,h,\gamma$ are uniquely determined by the set.
    
    First, we notice that $\gamma$ is the ``bottom-left'' corner of $\Q^{k,h}_\gamma$, since it is obtained by translating the dilation of a fundamental Riemannian child of $T$, and $A_j(0)=0$ for all $j=1,\ldots,\dimZ$. Hence, 
    \[
    \gamma_{j}:=\min_{\lambda \in \Q^{k,h}_\gamma}\lambda_j, \quad j=1,\ldots,\dimV,
    \] 
    and
    \[
    \gamma_j:=\min_{\lambda \in \Q^{k,h} \cap_{i=1,...,\dimV} \graffe{\lambda_i=\gamma_i} }\lambda_j, \quad j=\dimV + 1,\ldots,\dimV+\dimZ.
    \] 
    Then, if we fix the coordinates $x_2,\dots,x_{\dimV}$ and $z_{1}\dots,z_{\dimZ}$ of the points in $\Q^{k,h}_\gamma$, the first coordinate of the horizontal layer $x_1$ ranges over an interval of length $2^{-k-h}$, since $\Q^{k,h}_\gamma$ is obtained by dilating the $x$-components of $T$ by a factor of $2^{-k}$, and then splitting $h$ times by $2$ to obtain the $h$-th Riemannian generation. Hence, the number $k+h$ is uniquely determined. On the other hand, if we fix the coordinates $x_1,\dots,x_{\dimV}$ and $z_{2}\dots,z_{\dimZ}$ of the points in $\Q^{k,h}_\gamma$, the first coordinate of the vertical layer $z_1$ ranges over an interval of length $4^{-k} 2^{-h-1}$, since we split the vertical sides of $T$ $h$ times by $2$ when determining the Riemannian generation, and the homogeneous dilation contributes a factor of $4^{-k}$. This yields that $2k+h$ is also uniquely determined, and hence $h$ and $k$.
\end{proof}

Finally, we discuss the metric properties of the Riemannian tiles. Recall that the Riemannian distance that we consider on $N$ is also left-invariant, as is the sub-Riemannian one, but it is not homogeneous with respect to the dilations $\dil_t$. However, the Riemannian distance is equivalent to the sub-Riemannian distance at large scales, while at small scales it behaves like the Euclidean one.

\begin{proposition}
\label{prop-riemannian-metric}
For every compact subset $I$ of $\RR^+$ and every $t\in I$, $h \in \NN$ and $\gamma \in N$, 
there exists a point $n_\gamma^{\log_4t,h} \in \gamma \Nprod \dil_{t}(T^{\mathrm R,h})$ and constants $0<c_{\mathrm{R,I}}\leq C_{\mathrm{R,I}} <\infty $ such that 
\begin{equation} \label{eq-both}
    \riemball{n_\gamma^{\log_4t,h}}{c_{\mathrm{R,I}}  2^{-h}} \subseteq\gamma\Nprod \dil_{t}(T^{\mathrm R,h}) \subseteq \riemball{n_\gamma^{\log_4t,h}}{C_{\mathrm{R,I}} 2^{-h}}
\end{equation} where $c_{\mathrm{R,I}},C_{\mathrm{R,I}}$ are independent of $h$, $\gamma$ and $t$, but depend on $I$.
\end{proposition}
\begin{remark}
    If in the proposition above $t=4^{-k}$ for some $k \in \ZZ,$ then $\gamma\Nprod \dil_{t}(T^{\mathrm R,h}) = \Q^{k,h}_\gamma.$
\end{remark}
\begin{proof}
We assume $\gamma = (0,0)$, since the general case comes from the left invariance of the Riemannian metric. Fix $h\in \NN$ and $t\in I$. Recall that
\begin{equation*}
T^{\mathrm R,h} = \left\{ (x,z) \in N  :  2^h x \in [0,1)^{\dimV}, 
4^h (2z_j) \in \left[ A_j(2^h x),  A_j(2^h x) + 2^h\right), j = 1,\dots,\dimZ\right\},
\end{equation*} 
and hence 
\begin{align*}
\dil_t\tonde{T^{\mathrm R,h}} = \Big\{& (x,z) \in N \ : \ \frac{2^hx}{t} \in [0,1)^{\dimV}, \\ 
&\frac{4^h (2z_j)}{t^2} \in \left[ A_j(2^h x),  A_j(2^h x) + 2^h\right), \ j = 1,\dots,\dimZ \Big\}.
\end{align*} 
Let $A := \sup_j \|A_j\|_{L^{\infty}([0,1)^{\dimV})} < \infty$. If $h > \log_2(2A)$, 
we then have the two-sided inclusion
\begin{equation}
\begin{aligned}
\label{riemannian children}
[0,t\cdot 2^{-h})^{\dimV} \times \left[\tfrac {t^2} 2A \cdot 4^{-h}, -\tfrac {t^2} 2A \cdot 4^{-h} + t^2 \cdot 2^{-h-1}\right)^{\dimZ} \subseteq \dil_t(T^{\mathrm{R},h}) \\
\subseteq [0,t\cdot 2^{-h})^{\dimV} \times \left[-\tfrac {t^2} 2A \cdot 4^{-h}, \tfrac {t^2} 2A \cdot 4^{-h} + t^2 2^{-h-1}\right)^{\dimZ}.
\end{aligned}
\end{equation}
This happens whenever the left hand side interval in the $z$-direction is non-empty, which requires \( 2A \cdot 4^{-h} < 2^{-h},\) that is, $h > \log_2(2A)$. 

Therefore, both the inner and outer boxes in \eqref{riemannian children} are comparable to the Euclidean cubes whose side lengths is of order $2^{-h}$ in both components. Since the Riemannian distance on $N$ 
is equivalent to the Euclidean $\norma{\argomento}_\infty$ on $\RR^{\dimV+\dimZ}$ on small scales, and $t$ is bounded and $h\in \NN$, there exists a centre point $n_0^{\log_4t,h} \in \dil_t(T^{\mathrm{R},h})$ and constants $c_{\mathrm{R,I}}, C_{\mathrm{R,I}} > 0$, independent of $h$ but dependent on the interval $I$, such that \eqref{eq-both} holds when $h$ is sufficiently big.

On the other hand, if $0\leq h \leq \log_2(2A)$ (i.e. finitely many  values of $h$), \eqref{eq-both} is trivially true since $\dil_t(T^{\mathrm R,h}) $ are bounded and have non-empty interiors (this is a consequence of \cite[Corollary 2.3]{S92} and Proposition \ref{prop-riemannianchristcubes}).

The general case $\gamma \Nprod \dil_t(T^{\mathrm R,h})$ follows with $n_\gamma^{\log_4t,h} := \gamma \Nprod n_0^{\log_4t,h}$. 
\end{proof}

The fact that $(N, \riemdist)$ is doubling and Proposition \ref{prop-riemannian-metric} yield the following fact for the \emph{$\costante$-enlarged cube}
$$\left(\Q^{0,h}_\gamma\right)^* \defeq\left\{ n \in N \ : \ \riemdist(n, \Q^{0,h}_\gamma) <\costante \cdot \diam_{\mathrm{R}}(\Q^{0,h}_\gamma) \right\}.$$

\begin{corollary} 
Let $h\in \NN$, $\gamma \in N$ and $\costante >0$ a constant, then
 there exists a constant $\costantee>0$ independent of $h$ such that $$ \mu\tonde{\tonde{\Q^{0,h}_\gamma}^*} \leq \costantee \cdot \mu \tonde{\Q^{0,h}_\gamma}.$$
\end{corollary}
\begin{remark}
We can formulate a statement completely analogous to Remark \ref{rmk-christ} for the construction of Hyt\"onen and Kairema in the setting of Riemannian cubes endowed with the Riemannian metric. Since the Riemannian metric is equivalent to the Euclidean metric at small scales and to the sub-Riemannian metric at large scales, we fix a generation \(k\in\ZZ\) (for convenience, say \(k=0\)), corresponding to the fundamental tile. We then consider the nested partitions induced by the corresponding cubes at generations \((-hM_1,0)\), for \(h\in\NN\), and \((0,-hM_2)\), for \(h\in\NN\), where \(M_1\) and \(M_2\) are the integers arising from condition (2.3) in \cite[Theorem 2.2]{HK12}, so that some generations are skipped.

Although some technicalities may arise from the equivalence constants of the metrics involved, the argument proceeds exactly as in Remark \ref{rmk-christ}. Indeed, by applying to Theorem \ref{Thm-RiemannianChristcubes} and Proposition \ref{prop-riemannian-metric}, the proof carries over without further changes.
\end{remark}

\section{Admissible sets on Damek--Ricci spaces}\label{Sect-adminS}

We now construct the admissible sets in the CZ theory on $S$. Since these spaces are the semidirect product of an $H$-type group $N$ and $\RR$, it is natural to define admissible sets as products of two suitable sets, one in $N$ and one in $\RR$, subject to a compatibility condition. This condition can be interpreted as follows: when we translate such sets near the identity $e$, the logarithm of the size of the cube in $N$ and the length of the interval in $\RR$ must be comparable. The condition will be made precise in subsequent definitions.

Let $I$ be an interval of $\RR$ closed on the left and open on the right, and let $u_I$ be its midpoint, so that $I = [u_I - \rI, u_I + \rI)$, where $2\rI$ is the length of $I$. We say that $\widetilde{I} \in \mathcal{C}_\RR(I)$ if $\widetilde{I} = [u_I - \rI, u_I)$ or $\widetilde{I} = [u_I, u_I + \rI)$, which are the two disjoint intervals of equal length that partition $I$. We can repeat this procedure, obtaining the $\ell$-th generation $\mathcal C ^\ell_{\RR}(I)$ of $I$ by splitting it into two equal parts $\ell$ times.

For the subsequent discussion, we define the following:
\begin{itemize}
\item for \(\ell \in \ZZ\)
\begin{equation}\label{Def-mathcalIell}
	\mathcal{I}^{\ell}:=\{[u_{I}-2^{-(\ell+1)}, u_I+2^{-(\ell+1)})\  : \ u_I\in \RR\}.
\end{equation}
For convenience of notation, we will often denote, for $I \in \mathcal I ^\ell,$ $\rI:=2^{-(\ell+1)}$;

\item for $k\in\ZZ$ 
	\begin{equation}\label{Def-mathcalDk}
	\mathcal{D}^{k}:=\{\Q^k_\gamma=\gamma\Nprod \dil_{4^{-k}}(T) \ : \ \gamma\in N\};
	\end{equation}
    \item  for $k\in \ZZ$ and $h\in \NN,$ 
$$
\mathcal D^{k,h}:=\graffe{\Q^{k,h}_{\gamma}=\gamma \Nprod \diltile{-k}{T^{\mathrm{R},h}} \ : \ \gamma \in N}.
$$
\end{itemize}
\begin{definition}\label{Def-largeadmissiblesets}
Let $\mathcal{\Q} \in \mathcal{D}^k$ and $I \in \mathcal{I}^{\ell}$, for $k \in \ZZ$ and $\ell \leq 0$ (so $\rI := 2^{-(\ell +1)} \geq \tfrac{1}{2}$).
We call a set of the form $\mathcal{\Q} \times I=\Q^k_\gamma \times [u_I-\rI, u_I+\rI) \subseteq S$ a \emph{large admissible set} if it satisfies the following admissibility condition \begin{equation}\label{Ieq-admissiblecondition}
\ee^{\frac{u_I}{2}}\ee^{\rI}\leq2^{-k}<\ee^{\frac{u_I}{2}}\ee^{\frac92\rI}.
\end{equation}
\end{definition}
Now we want to define the children of large admissible sets.
The idea is the following: we can split a large admissible set $\mathcal{\Q} \times I$ in two different ways. The first way consists of considering the cubes $\widetilde{\mathcal{\Q}} \times I$, given by the product of the $2^{\dimV + 2\dimZ}$ sub-Riemannian children of $\widetilde{\Q} \in \mathcal{C}_{\mathrm{SR}}(\Q)$ with the interval $I$, provided the $\widetilde{\Q}\times I$
satisfy condition \eqref{Ieq-admissiblecondition} (notice that if one of the $\widetilde{\Q} \in \mathcal{C}_{\mathrm{SR}}(\Q)$ satisfy the condition, then all the others do it, since the condition only depends on the size of the cubes and children of the same cube share the same size). 

When this condition fails, we resort to the second way of splitting $\mathcal{\Q} \times I$, namely by taking the two children of the interval $I$, $\widetilde{I} \in \mathcal{C}_\RR(I)$ and forming their product $\Q \times \widetilde{I}$. 

Eventually, one reaches the limiting scale \(I\in\mathcal I^0\), at which \(I\) can no longer be subdivided while remaining in the large-scale regime. If, at the same time, \(\widetilde{\Q}\times I\) fails to satisfy \eqref{Ieq-admissiblecondition} for every sub-Riemannian child \(\widetilde{\Q}\) of \(\Q\), neither of the two preceding splitting procedures produces large admissible sets. We call such a set a \emph{transition admissible set}; its precise definition is given below.

Before giving this definition, we establish a lemma showing that as long as $\ell \leq -1$, i.e., the length of the interval $I$ is greater than or equal to $2$, the set $\mathcal{\Q} \times I$ can always be split in one of the two ways described above, and the children are still large admissible sets.

\begin{lemma}\label{Lem-childrenexistence}
Let \(\Q\times I = \Q^k_\gamma \times [u_I-\rI, u_I +\rI)\) be a large admissible set with \(\ell\le -1\) (equivalently \(\rI =2^{-(\ell+1)}\ge 1\)).  
If for every \(\widetilde \Q\in\mathcal{C}_{\mathrm{SR}}(\Q)\) the product \(\widetilde \Q\times I\) is \emph{not} a large admissible set, then for every \(\widetilde I\in\mathcal{C}_\RR(I)\) the product \(\Q\times\widetilde I\) is a large admissible set.
\end{lemma}

\begin{proof}
Since \(\Q\times I= \Q^k_\gamma \times [u_I-\rI, u_I +\rI)\) is large admissible set, we have  
\begin{equation}\label{Inq-adm1}
    \ee^{\frac{u_I}{2}}\ee^{\rI}\;\le\;2^{-k}\;<\;\ee^{\frac{u_I}2}\ee^{\frac{9}{2}\rI}.
\end{equation}
For any \(\widetilde \Q\in\mathcal{C}_{\mathrm{SR}}(\Q)\), we have \(\widetilde \Q\in\mathcal{D}^{k+1}\); thus, \(\widetilde{\Q}\times I\) would be a large admissible set exactly when \(\ee^{\frac{u_I}{2}}\ee^{\rI}\le 2^{-(k+1)}\).  
By hypothesis this fails for every \(\widetilde \Q\), hence  
\begin{equation}\label{Inq-adm2}
2^{-(k+1)} < \ee^{\frac{u_I}{2}}\ee^{\rI}.
\end{equation}
Now take any $\widetilde{I} \in \mathcal{C}(I)$. Then $\widetilde{I}$ has midpoint $u_{\widetilde{I}} = u_I \pm \frac{\rI}{2}$ and half length $\rI' = \frac{\rI}{2}$. We need to verify that
\[
\ee^{\frac{u_{\widetilde{I}}}{2}}\ee^{\rI'} \leq 2^{-k} < \ee^{\frac{u_{\widetilde{I}}}{2}}\ee^{\frac{9}{2}\rI'},
\]
which would mean that both children $\mathcal{\Q} \times \widetilde{I}$ satisfy the admissibility condition \eqref{Ieq-admissiblecondition}.

As for the lower bound, we have
\[
\ee^{\frac{u_{\widetilde I}}{2}}\ee^{\rI'}=\ee^{\frac{u_{I}}{2}}\ee^{\pm \frac{\rI}{4}}\ee^{\frac{\rI}{2}}
\le \ee^{\frac{u_I}{2}}\ee^{\rI}\le 2^{-k},
\]
where we used \eqref{Inq-adm1}, so the lower bound holds. For the upper bound, we compute
\[
\ee^{\frac{u_{\widetilde{I}}}{2}}\ee^{\frac{9}{2}\rI'} = \ee^{\frac{u_I}{2}}\ee^{\frac{9\rI}{4} \pm \frac{\rI}{4}}\ge \ee^{\frac{u_I}{2}}\ee^{2\rI}.
\]
Hence, we need to verify that
\[
2^{-k} < \ee^{\frac{u_I}{2}}\ee^{2\rI}.
\]
Since $\rI \geq 1$, we have $\ee^{2\rI} \geq 2\ee^{\rI}$ (because $\ee^{\rI} \geq \ee > 2$). Using \eqref{Inq-adm2},
\[
2^{-k} < 2\,\ee^{\frac{u_I}{2}}\ee^{\rI} \leq \ee^{\frac{u_I}{2}}\ee^{2\rI},
\]
so the required inequality holds.

Therefore \(\Q\times\widetilde I\) satisfies the admissibility condition \eqref{Ieq-admissiblecondition} for both the children of $I$. Moreover, \(\ell+1\le0\) (i.e. \(\rI'=\frac{\rI}{2}\ge \frac{1}{2}\)), so it is indeed a large admissible set.
\end{proof}

The iterative process, together with Lemma \ref{Lem-childrenexistence}, leads to the following notion of a transition set, which occurs when the sub-Riemannian children of $\mathcal{\Q}$ fail to satisfy the admissibility condition \eqref{Ieq-admissiblecondition} and, at the same time, the length of $I$ equals exactly $1$.
\begin{definition}\label{Def-transitionsets}
		A \emph{transition admissible set} is a set of the form $\Q\times I=\Q^k_\gamma \times [u_I-\rI, u_I+\rI)\in \mathcal{D}^{k}\times \mathcal{I}^0$ that satisfies the following condition \begin{equation}\label{Ieq-transitioncondition}
    2^{-k-1}<\ee^{\frac{u_I}{2}}\ee^{\frac12}\leq2^{-k}.
    \end{equation}
	\end{definition}
Notice that transition admissible sets are large admissible sets. For these sets, we must find a new way of splitting into children: namely, we shall split dyadically simultaneously $\Q$ and $I$, where in the case of $\Q$ we shall use Riemannian descendants. We can turn to the definition of small admissible sets.
\begin{definition}\label{Def-smalladmissiblesets} 
A \emph{small admissible set} is a set of the form $\Q \times I \subseteq S$, for some $\Q \in \mathcal C^h_{\mathrm R} (\widehat{\Q})$ and $I \in \mathcal C^h_{\RR} (\widehat{I})$, where $\widehat{\Q} \times \widehat{I} $ is a \emph{transition admissible set} for some \(h\in\NN\).
\end{definition}
\begin{remark}
We observe that the definition above can be also given in a quantitative way as follows: a small admissible set is a set of the form $\Q \times I \subseteq S$, with $\mathcal{\Q} \in \mathcal{D}^{k,h}$, $I \in \mathcal{I}^{h}$, for $k \in \ZZ$ and $h \in \mathbb{N}$, and such that there exists 
\begin{equation}
\label{eq-smallcondition}
u_{\widetilde{I}} \in \left\{ u_I + 2^{-h-1} - \tfrac{1}{2} + \ell\, 2^{-h} \right\}_{\ell = 0, \ldots, 2^h - 1},
\end{equation}
satisfying the admissibility condition
\[
2^{-k-1} < \mathrm{e}^{\frac{u_{\widetilde{I}}}{2}} \mathrm{e}^ {\frac 1 2 } \leq 2^{-k}.
\]
\end{remark}
Children of small admissible set $\Q\times I$ are declared to be $\widetilde{\Q}\times \widetilde{I}$, where $\widetilde{\Q}\in\mathcal{C}_{\mathrm{R}}(\Q)$ and $\widetilde{I}\in \mathcal{C}_\RR(I)$. The admissibility condition is still satisfied: since we split simultaneously $Q$ and $I$, we have that $\widetilde{\Q} \in \mathcal{C}^{h+1}_{\mathrm{R}}(\widehat \Q)$ and $\widetilde{I}\in \mathcal{C}^{h+1}_\RR(\widehat{I})$ in Definition \ref{Def-smalladmissiblesets}.

\begin{definition}\label{Def- S generation}
Let $R\subseteq S$ be an admissible set (either large or small). We say that $R_j\in \graffe{R_i}_{i=1}^{n}$ belongs to the \textit{first generation} of $R$, or, equivalently, is one of the \emph{children} of $R$, and write
$R_j \in\mathcal{C}(R)= \mathcal{C}^1(R)$, if
$R_j$ is obtained from $R$ through the splitting process described above.
More generally, for $\ell\geq 2$, we define the \emph{$\ell$-generation} of $R$ recursively as
$$\mathcal C^\ell(R)= \mathcal C^1(\mathcal C^{\ell-1}(R)).$$
We also define $$\mathcal C^{\ast} (R) := \bigcup_{k\in\NN} \mathcal C^k(R).$$
\end{definition}
\begin{remark}
    The number of children of an admissible set is either $2^{\dimV+2\dimZ}$, $2$, (for large admissible sets which are not transition admissible sets) or $2^{\dimV+\dimZ+1}$ (for transition and small admissible sets).
\end{remark}
\begin{proposition}
\label{prop-regularity}
    Let $R=\Q^{k,h}_\gamma \times[u_I-\rI, u_I +\rI)$ be an admissible set and let $R_j\in \mathcal C (R).$ Then
    \begin{equation}
    \label{eq-number of children}
        \frac{\rho(R)}{\rho(R_j)}=\begin{cases}
            2^{\dimV+\dimZ+1} \quad &\mathrm{if \ } R \mathrm{\ is\ small \ or \ transition\ admissible \ set,}\\
            2^{\dimV+2\dimZ}  &\mathrm{if \ } 2^{-k-1}\geq \mathrm{e}^{\frac{u_I}{2}}\mathrm{e}^{\rI}\ \mathrm{and \ } h=0,\\
            2 &\mathrm{otherwise.  }
        \end{cases}
    \end{equation}
\end{proposition}
\begin{proof} Since the children of \(R\) form a partition of R into sets of equal \(\rho\)-measure (from Fubini's Theorem), the ratio \(\frac{\rho(R)}{\rho(R_j)}\) equals to the number of children of \(R\). Conclusion \eqref{eq-number of children} follows from the way of splitting admissible sets.
\end{proof}
    \begin{remark}
    \label{rmk-splitting}
    If we start splitting a small admissible set, we obtain small admissible sets. On the other hand, if we start from a large admissible set, we eventually reach small admissible sets as well. Indeed, a large admissible set can be split in the $N$-direction (increasing $k$ by $1$ at each step) only finitely many times, since after finitely many such steps the left-hand inequality in the admissibility condition fails. At that point, the set must instead be split in the $\mathbb R$-direction, and this process is repeated finitely many times, until we reach a transition admissible set whose children are small; from then on, all of its descendants are small as well.
\end{remark}

In order to investigate some geometric properties of the above defined sets we need a few technical results.
\begin{lemma}\label{Lem-inclusions}
    For any $\mathbf{x}=(x,z,u)\in S$ and $r\geq 1$, the following hold.
	\begin{enumerate}
	\item \label{lem4.11-item1} There exists \(\kappa_1>0\), such that 
	$$
	  \damekricciball{e}{r}\subseteq \subriemball{(0,0)}{\kappa_1\ee^{\frac{r}{2}}}\times [-r,r].
	$$
	\item \label{lem4.11-item2} For the same \(C_N\) as in \eqref{Equ_inclusion}, there exist a positive constant \(\kappa_2\) such that 
	$$
	  \subriemball{(0,0)}{C_N\mathrm{e}^{\frac{9}{2}r}}\times[-r,r)\subseteq \damekricciball{e}{\kappa_2r}.
	$$
	\end{enumerate}
\end{lemma}
\begin{proof}
	\eqref{lem4.11-item1} From Lemma \ref{Lem-coordinateestimates}, we get
	\begin{equation}\label{Inclu-ax+b}
	  \damekricciball{e}{r}\subseteq \horiball{0}{2\sinh\frac{r}{2}}\times \vertball{0}{\sinh r}\times[-r,r].
	\end{equation}
    Since 
    \[
    \sinh r=2\sinh\frac{r}{2}\cosh \frac{r}{2} =2\coth\frac{r}{2} \sinh^2 \frac{r}{2} \leq 2\coth\frac{1}{2}\sinh^2 \frac{r}{2}
    \]
    for \(r\geq 1\), then \eqref{Eq-ball-box-N} implies that, there exists a constant \(\kappa_1>0\) such that 
    \[
   \horiball{0}{2\sinh\frac r2}\times  \vertball{0}{\sinh r}\subseteq\subriemball{(0,0)}{\kappa_1\sinh\frac r2}.
   \]
   Together with \eqref{Inclu-ax+b}, it concludes point \eqref{lem4.11-item1}.\\
	\eqref{lem4.11-item2} By \eqref{Eq-distanceinS} and the bijection formula for \(\cosh\), 
	$$
	\begin{aligned}
		\cosh |\mathbf{x}|&=2\tonde{\cosh\frac{u}{2}+\frac{1}{2}\mathrm{e}^{-\frac{u}{2}}\left|\frac{x}{2}\right|^2}^2+\frac{1}{2}\mathrm{e}^{-u}|z|^{2}-1\\
		&\leq4\left[\cosh^2\frac{u}{2}+\frac{1}{4}\mathrm{e}^{-u}\left|\frac{x}{2}\right|^4\right]+\frac{1}{2}\mathrm{e}^{-u}|z|^{2}-1\\
		&=4\cosh^2\frac{u}{2}+\frac{1}{16}\mathrm{e}^{-u}|x|^{4}+\frac{1}{2}\mathrm{e}^{-u}|z|^{2}-1\\
		&=2\cosh u+1+\frac{1}{16}\mathrm{e}^{-u}|x|^{4}+\frac{1}{2}\mathrm{e}^{-u}|z|^{2}.
	\end{aligned}
	$$	
	Since $r\geq 1$, there exist \(\kappa_2\) such that 
	\begin{align*}
		2\cosh r+1&\leq \frac{1}{3}\cosh (\kappa_2r), \\ 
		\frac{1}{16}\mathrm{e}^{ r}|\beta_1C_{N}\mathrm{e}^{\frac{9}{2}r}|^4&\leq\frac{1}{3}\cosh (\kappa_2r),\,\, \\ 
		\frac{1}{2}\mathrm{e}^{r}|\beta_1\,C^2_{N}\mathrm{e}^{9r}|^2&\leq\frac{1}{3}\cosh (\kappa_2r),
	\end{align*}
    where \(C_N\), \(\beta_1\) are as in \eqref{Eq-ball-box-N}.
	Thus, we get from \eqref{Eq-ball-box-N} that if $\mathbf{x}\in \subriemball{(0,0)}{C_{N}\mathrm{e}^{\frac{9}{2}r}}\times[-r,r]$, then 
	$$
	\cosh |\mathbf{x}|\leq \cosh (\kappa_2 r),
	$$
	which implies that $\mathbf{x}\in \damekricciball{e}{\kappa_2r}$.
\end{proof}

For any admissible set \(R\), let
\begin{equation*}
R^\ast\defeq\left\{\mathbf x\in S:\damekriccidist(\mathbf x,R)\leq\frac{1}{2\kappa_2}\cdot\diam_{\mathrm S}(R)\right\},
\end{equation*}
be the set defined in \eqref{Def-enlargedset} with \(\constantS=\frac{1}{2\kappa_2}\), and \(\kappa_2\) from Lemma \ref{Lem-inclusions}. In the following proposition we will show that \(\rho(R)\) is comparable to \(\rho(R^*)\) for any large admissible set \(R\), while the case of small admissible sets will be treated in Proposition \ref{Propo-smallinclusions}.
\begin{proposition}\label{Propo-largeinclusions}
	Denote by $R=\Q\times I=\Q_{\gamma}^{k}\times[u_I-\rI, u_I+\rI)$ a large admissible set in $\mathcal{D}^{k}\times \mathcal{I}^{\ell}$. Let $c_N$ and $C_N$ be as in Theorem \ref{Thm-Christcubes} and \(\kappa_1\), \(\kappa_2\) in Lemma \ref{Lem-inclusions}. Then the following hold:
	\begin{enumerate}
		\item \label{prop-item1}$R\subseteq \damekricciball{(n_{\gamma}^{k},u_{I})}{\kappa_2\rI}$, where $n^k_\gamma$ is as in \eqref{teo3.5-item4} of Theorem \ref{Thm-Christcubes};
	
		\item \label{prop-item2}\(\rho(R)\) is comparable with \(\rho(R^*)\), more precisely
		$$
		\rho(R^*)\leq 2\left(\frac{\kappa_1+C_N}{c_N}\right)^{2\dimHomog}\rho(R).
		$$
	\end{enumerate}
\end{proposition}
\begin{proof}
    \eqref{prop-item1} By \eqref{Equ_inclusion}, we know that 
    $$
    \begin{aligned}
        R&\subseteq \subriemball{n_\gamma^k}{C_N2^{-k}}\times [u_I-\rI, u_I+\rI)\\
        &\subseteq (n_\gamma^k,u_{I})\Sprod \left( \subriemball{(0,0)}{C_{N}\ee^{-\frac{u_{I}}{2}}2^{-k}}\times [-\rI,\rI)\right),
    \end{aligned}
    $$
 so, from the admissibility condition \eqref{Ieq-admissiblecondition} and \eqref{lem4.11-item2} in Lemma \ref{Lem-inclusions},  
 $$
 \begin{aligned}
     R&\subseteq(n_\gamma^k,u_{I})\Sprod \left( \subriemball{(0,0)}{C_{N}\mathrm{e}^{\frac{9\rI}{2}}}\times [-\rI,\rI)\right)\\
     &\subseteq (n_{\gamma}^{k},u_{I})\Sprod \damekricciball{e}{\kappa_2\rI}=\damekricciball{(n_{\gamma}^k,u_I)}{\kappa_2\rI}.
 \end{aligned}
 $$ 
	
	\eqref{prop-item2} Notice that from \eqref{prop-item1}, it follows that \(\diam_{\mathrm{S}}(R)\le 2\kappa_2 \rI\le \rI/\constantS\). So, by the definition of the enlargement and left translation, we know that 
	$$
	\begin{aligned}
		R^{*}:&=\{\mathbf{x}\in S:\damekriccidist(\mathbf{x},R)<\constantS\cdot\diam_\mathrm{S}(R)\}
		\subseteq\bigcup_{\mathbf{x}\in R}\damekricciball{\mathbf{x}}{\rI}=\bigcup_{\mathbf{x}\in R}\mathbf{x}\Sprod \damekricciball{e}{\rI}.
	\end{aligned} 
    $$
     By \eqref{lem4.11-item1} in Lemma \ref{Lem-inclusions} and
     \eqref{Equ_inclusion}, we have that 
    $$
    \begin{aligned}
    	R^*&\subseteq\bigcup_{\mathbf{x}=(x,z,u)\in R}\{(x,z,u)\Sprod\left(\subriemball{(0,0)}{\kappa_1\,\ee^{\frac{\rI}{2}}}\times[-\rI,\rI]\right)\}\\
    	&= \bigcup_{(x,z,u)\in R}\left(\subriemball{(x,z)}{\kappa_1\,\ee^{\frac{u}{2}}\ee^{\frac{\rI}{2}}}\times[u- \rI,u+\rI]\right)\\
    	&\subseteq \subriemball{n_\gamma^k}{\kappa_1\,\ee^{\frac{u_I+\rI}{2}}\ee^{\frac{\rI}{2}}+C_N2^{-k}}\times[u_{I}-2\rI,u_{I}+2\rI]\\
    	&\subseteq \subriemball{n_\gamma^k}{(\kappa_1+C_N)2^{-k}}\times[u_{I}-2\rI,u_{I}+2\rI]),
    	\end{aligned}
    $$
    where we used that 
    \[
    R\subseteq \subriemball{n_\gamma^k}{C_N2^{-k}}\times [u_I-\rI, u_I+\rI]
    \]
    and the last inclusion comes from the admissibility condition \eqref{Ieq-admissiblecondition}. Applying \eqref{Equ_inclusion} again, we deduce that 
    $$
    \rho(R^*)\leq \left(\frac{\kappa_1+C_N}{c_N}\right)^{2\dimHomog}|\subriemball{n_\gamma^k}{c_N2^{-k}}|\cdot 4\rI\leq 2\left(\frac{\kappa_1+C_N}{c_N}\right)^{2\dimHomog}\rho(R). 
    $$
\end{proof}
The following lemma, together with the subsequent remark, will be used both to construct an exhausting sequence of admissible sets in Theorem \ref{Thm-CZspaceproof} and later in Section \ref{Sec-endpointRiesz}.
\begin{lemma}\label{Lem-compactinadmissible}
    For every \(r>0\), there exists a large admissible set \(R\in\admis\) such that \(\riemball{e}{r}\subseteq R\).
\end{lemma}
\begin{proof}
    Without loss of generality, we can prove only the case \(r\geq 1\). By Lemma \ref{Lem-inclusions} \eqref{lem4.11-item1},
    \[
        \damekricciball{e}{r}\subseteq \subriemball{(0,0)}{\kappa_1\ee^{r/2}}\times[-r,r).
    \]
Choose \(\ell\leq0\) sufficiently small so that, letting \(\rI=2^{-(\ell+1)}\), 
and for \(c_N\) in \eqref{Equ_inclusion}, we have
    \begin{equation}\label{Ieq-compact1}
        \rI\geq r \,\,\,\text{and}\,\,\,
        c_N\ee^{\rI}\geq\kappa_1\ee^{\frac{r}{2}}.
    \end{equation}
    Set \( I=[-\rI,\rI)\): \(I\in \mathcal{I}^{\ell}\) (by \eqref{Def-mathcalIell}) and \(u_I=0\). Next, choose \(k\in\ZZ\) such that
    \begin{equation}\label{Ieq-compact2}
        \ee^{\frac{u_I}{2}}\ee^{\rI}=\ee^{\rI}\leq2^{-k}<\ee^{\frac{9\rI}{2}}=\ee^{\frac{u_I}{2}}\ee^{\frac{9\rI}{2}}.
    \end{equation}
    For \(n_0^0 \) in Proposition \ref{subriemannian balls}, define
    \[
        \gamma\defeq\bigl(\dil_{4^{-k}}(n_0^0)\bigr)^{-1},\,\,\,
    \Q_\gamma^k\defeq \gamma\cdot_{N}\dil_{4^{-k}}(T).\]
    According to the definition of \(\mathcal D^k\) in \eqref{Def-mathcalDk}, \(\Q_\gamma^k\in \mathcal D^k\), and the centre of \(\Q_\gamma^k\) is \(\gamma\cdot_N n^0_0=(0,0)\). Therefore, by \eqref{Equ_inclusion}, we deduce that 
    \[
        \subriemball{(0,0)}{c_N2^{-k}}\subseteq Q_\gamma^k.
    \]
    Moreover, combining \eqref{Ieq-compact1} and \eqref{Ieq-compact2},
    \[
        c_N2^{-k}\geq c_N\ee^{\rI}\geq\kappa_1\ee^{\frac r2}.
    \]
    Consequently,
    \[
        \damekricciball{e}{r}\subseteq \Q_\gamma^k\times I.
    \]
    \eqref{Ieq-compact2} yields that \(\Q_\gamma^k\times I\) is a large admissible set, according to Definition \ref{Def-largeadmissiblesets}.
\end{proof}
\begin{corollary}\label{Rmk-compactinadmissible}
    For every compact set \(K\subseteq S\), 
    there exists a large admissible set \(R\in\admis\) such that \(K\subseteq R\).
\end{corollary}

Let us now turn to small admissible sets. We want to outline some properties analogous to those stated above for large admissible sets. 
\begin{proposition} \label{Propo-smallinclusions}
    	Denote by $R=\Q\times I=\Q_{\gamma}^{k,h}\times[u_I-\rI, u_I+\rI)$ a small admissible set in $\mathcal{D}^{k,h}\times \mathcal{I}^{h}$. For some positive constants $C_{\mathrm{S}}$ and $C_\infty$, depending only on $S$, the following hold: 
	\begin{enumerate}
		\item \label{prop4.15-item1} there exists $\mathbf{x}_R \in R$ such that 
        \begin{equation}
        \label{eq-palle ammissibili}
             R \subseteq B_S((\gamma,u_I), C_{\mathrm{S}} \rI);
        \end{equation}
		\item \label{prop4.15-item2}for the measure of $R$, we have the estimate: 
        $$
        2^{\dimV+1} \ee^{\frac{\dimHomog}{2}} \ee^{\dimHomog u_I} \rI^{\dimV+\dimZ+1}<\rho(\Q^{k,h}_\gamma\times I)< 2^{2(\dimV+\dimZ)+1} \ee^{\frac{3\dimHomog}{2}} \ee^{\dimHomog u_I} \rI^{\dimV+\dimZ+1};
        $$
		\item \label{prop4.15-item3}we have the estimate 
		$$ \rho(R^* ) \leq  C\,\tonde{(1+4\constantS)C_\infty C_{\mathrm{S}}}^{\dimV+\dimZ+1}2^{-1-\dimV}\ee^{-\frac{\dimHomog}{2}}\rho(R),$$
        for some constant \(C>0\) independent of \(h\) and \(k\).
	\end{enumerate}
\end{proposition}
\begin{proof}
    Let $R = \Q^{k,h}_\gamma \times I$ be a small admissible set as above, and hence $\rI=2^{-h-1}$. Then,
\begin{align}
   \notag R&= \tonde{\gamma\Nprod \diltile{-k}{T^{R,h}}}\times [u_I - 2^{-h-1},u_I+2^{-h-1}) \\
   \label{eq-translation}&=(\gamma,u_I)\Sprod \tonde{\dil_{\ee^{-u_I}4^{-k}}(T^{R,h})\times [ - 2^{-h-1},2^{-h-1})}.
    \end{align}
    Due to the admissibility condition for small admissible sets in Definition \ref{Def-smalladmissiblesets}, we have that 
    \[
    \ee \leq \ee^{-u_{\widehat{I}}} 4^{-k} < 4\ee,
    \]
    for some $h$-ancestor $\widehat{I}$ of $I$. Since $u_I \in I \subseteq \widehat I =\left[u_{\widehat I} - \tfrac 1 2, \ u_{\widehat I} + \tfrac 1 2\right) $, we have
    \begin{equation}\label{Inq-kdoubleestimates}
    \ee^{\frac{1}{2}}<  \ee^{-u_{I}}4^{-k}<4\ee^{\frac{3}{2}}.
    \end{equation}
    Therefore, we can see every small admissible set as 
    $$
    R=\textbf{x}_\mathrm{R}\Sprod\tonde{\dil_{t}(T^{R,h})\times[-2^{-h-1},2^{-h-1})}
    $$
    where $\textbf{x}_{\mathrm R}=(\gamma,u_I) \in S$ and \(t\in(\ee^{\frac{1}{2}},4\ee^{\frac{3}{2}})\).\\
\eqref{prop4.15-item1} Since the Damek--Ricci distance is left-invariant, it suffices to prove (1) just for
\begin{equation}
\label{eq-t non uno}
R_{t,h}:= \dil_t (T^{R,h}) \times \left[-2^{-h-1}, 2^{-h-1}\right),
\end{equation}
for \(t\in(\ee^{\frac{1}{2}},4\ee^{\frac{3}{2}})\).

For every compact set $K\subseteq S$, the Damek--Ricci distance and the
Euclidean distance are equivalent on $K$, by
\cite[Lemma 2.2]{CCM19}. Hence, after fixing a compact
set $K$ containing all the sets under consideration, there exist constants
$c_\infty,C_\infty>0$ and $r_0>0$ such that, for every
$\mathbf{x}\in K$ with $B_S(\mathbf{x},r_0)\subseteq K$ and every
$0<r\le r_0$, one has
\begin{equation}\label{Eq-equivalence-supdr}
     \supeuclball{\mathbf{x}}{c_\infty r}\subseteq
    \damekricciball{\mathbf{x}}{r}\subseteq  \supeuclball{\mathbf{x}}{C_\infty r}.
\end{equation}
Here
\[
    \supeuclball{\mathbf{x}}{r} \defeq
    \left\{
    \mathbf{x}\Sprod \mathbf{y}\in S:
    \|\mathbf{y}\|_{\infty}<r
    \right\},
\]
where, writing \(\mathbf{y}=(y_1,\dots,y_{\dimV+\dimZ+1})\),
\[
    \|\mathbf{y}\|_{\infty}\defeq
    \max_{1\le i\le \dimV+\dimZ+1}|y_i|.
\]
From \eqref{riemannian children} we deduce that
\begin{align*}
  R_{t,h}  \subseteq  
[0,t^{1/2}2^{-h})^{\dimV} \times \left[-2tA \cdot 4^{-h}, 2tA \cdot 4^{-h} + t\cdot 2^{-h-1}\right)^{\dimZ} \times\left[-2^{-h-1}, 2^{-h-1}\right),
\end{align*}
which is a rectangular box in $\RR^{\dimV+\dimZ+1}$ containing $(0,0,0)$. Since \(t\) ranges in a compact subinterval of \((0,\infty)\), namely
\(t\in(\ee^{\frac{1}{2}},4\ee^{\frac{3}{2}})\), both \(t^{\frac{1}{2}}\) and \(t\) are bounded above and below by positive constants independent of \(h\). Therefore, the side lengths of this outer box are all comparable to \(2^{-h}\), uniformly in \(h\) and $t$. 
Together with the local equivalence between the sup-norm boxes and the Damek--Ricci balls, this implies that there exists a constant \(C_{\mathrm{S}}>0\), independent of \(h\), such that
\begin{equation}\label{Inc-Rth}
R_{t,h}\subseteq\damekricciball{e}{C_{\mathrm{S}}2^{-h}}.
\end{equation}
\eqref{eq-palle ammissibili} comes from left translation by $\mathbf{x}_{\mathrm R} = (\gamma,u_I).$\\
\eqref{prop4.15-item2} 
From point \eqref{teo3.11-item2} in Theorem \ref{Thm-RiemannianChristcubes} and with \(\rI=2^{-h-1}\), we deduce that 
\[
\rho(\Q^{k,h}_\gamma\times I)=2^{-h(\dimV+\dimZ+1) - k(\dimV+2\dimZ)-\dimZ }.
\]
Using the estimates for \(k\) in \eqref{Inq-kdoubleestimates}, it follows that 
\begin{equation}\label{eq-bound on Rh}
2^{\dimV+1} \ee^{\frac{\dimHomog}{2}} \ee^{\dimHomog u_I} \rI^{\dimV+\dimZ+1}<\rho(\Q^{k,h}_\gamma\times I)< 2^{2(\dimV+\dimZ)+1} \ee^{\frac{3\dimHomog}{2}} \ee^{\dimHomog u_I} \rI^{\dimV+\dimZ+1}.
\end{equation}

\eqref{prop4.15-item3}  We prove it first for $R_{t,h}$ and $R^*_{t,h}$, and then extend it to every small admissible set. From \eqref{Inc-Rth}, we obtain \(\diam_{\mathrm{S}}(R_{t,h})\leq 4C_{\mathrm{S}} \rI \). Hence, combined with \eqref{Eq-equivalence-supdr} and the fact that \(\rI=2^{-h-1}\), we obtain
\[
R_{t,h}^*\subseteq\damekricciball{e}{(1+4\constantS)C_{\mathrm{S}} \rI}\subseteq \supeuclball{e}{(1+4\constantS)C_\infty C_{\mathrm{S}}\rI}.
\]

Therefore,
\begin{equation}\label{Est-rhoR^*}
\rho(R_{t,h}^*)
\leq
C\,\tonde{(1+4\constantS)C_\infty C_{\mathrm{S}}\rI}^{\dimV+\dimZ+1},
\end{equation}
for some constant \(C>0\) independent of \(t\) and \(h\). 

Similarly as \eqref{eq-bound on Rh}, we get
\[
\rho(R_{t,h})\geq 2^{-h(\dimV+\dimZ+1)-\dimZ} t^{\dimHomog}>\rI^{\dimV+\dimZ+1}2^{1+\dimV}\ee^{\frac{\dimHomog}{2}},
\]
since \(t>\ee^{\frac{1}{2}}\), which, together with \eqref{Est-rhoR^*}, yields
\begin{align*}
\rho(R_{t,h}^* ) &\leq C\,\tonde{(1+4\constantS)C_\infty C_{\mathrm{S}}\rI}^{\dimV+\dimZ+1}\tonde{\rI^{\dimV+\dimZ+1}2^{1+\dimV}\ee^{\frac{\dimHomog}{2}}}^{-1} \rho (R_{t,h})\\
&\leq C\,\tonde{(1+4\constantS)C_\infty C_{\mathrm{S}}}^{\dimV+\dimZ+1}2^{-1-\dimV}\ee^{-\frac{\dimHomog}{2}}\rho (R_{t,h})
\end{align*}
with the constant in the inequality independent of $h$ and \(t\).

Since the metric is left-invariant, if $R=\mathbf x \cdot_S R_{t,h}$, then
$R^*=\mathbf x \cdot_S R_{t,h}^*$. Moreover, left translation of sets by $\mathbf{x}$ multiplies the measure of such sets by the same factor (the modular function). That means,
$$
\frac{\rho(R^*)}{\rho(R)}=\frac{\rho(\mathbf{x}R^*)}{\rho(\mathbf{x}R)}, \quad \forall\mathbf{x} \in S,$$
and hence the inequality 
$$ \rho(R^* ) \leq  C\,\tonde{(1+4\constantS)C_\infty C_{\mathrm{S}}}^{\dimV+\dimZ+1}2^{-1-\dimV}\ee^{-\frac{\dimHomog}{2}}\rho(R)$$ still holds for every small admissible set $R$.
\end{proof}
\begin{corollary}
\label{cor-R contenuto in palla}
    For every admissible set $R=\Q\times I=\Q^{k,h}_\gamma\times[u_I-\rI, u_I+\rI)$ (large if $h=0$, small if $h>0$), we have 
    \begin{equation}
        \label{eq-palla}
        R\subseteq \damekricciball{\mathbf{x}_{\mathrm{R}}}{c \rI}
    \end{equation}
    for $c=\max\graffe{C_{\mathrm{S}},\kappa_2 }$ and $\mathbf{x}_{\mathrm{R}} = \begin{cases}
        (n^k_\gamma,u_I) \quad &\mathrm{if \ } h=0,\\
        (\gamma,u_I) &\mathrm{if \ } h>0.
    \end{cases}$
\end{corollary}
\begin{proposition}
\label{prop-discendenti piccoli}
    For all $\varepsilon>0$ and $R=\Q^{k,h}_\gamma\times[u_I-\rI, u_I+\rI)$ admissible set (either large or small), and for all $\mathbf x \in R$, there exists a descendant $R'=\Q^{k',h'}_{\gamma '}\times[u_I'-\rI', u_I'+\rI')$ of $R$ such that $\mathbf x \in R'$ and $\rI '\leq\varepsilon.$
\end{proposition}
\begin{proof}
Fix an admissible set \(R=\Q^{k,h}_\gamma\times[u_I-\rI,u_I+\rI)\), a point \(\mathbf x\in R\), and \(\varepsilon>0\). If
\(\rI\leq\varepsilon\), there is nothing to prove.

     We assume therefore that $\rI >\varepsilon$ and we build a chain of descendants of $R$, passing at each step to the unique child containing $\mathbf x$, and we track how $\rI$ evolves along the way. If $R$ is large and it splits along the $N$-side into finitely many admissible children, then exactly one of them contains $\mathbf x$, and we pass to it. However, this $N$-side splitting can only be performed finitely many times: since it increases $k$ at each step, and $2^{-k}$ decreases accordingly, the admissibility condition from below \eqref{Ieq-admissiblecondition} eventually fails, so after finitely many steps splitting along $N$ is no longer possible. At that point we switch to splitting along $\RR$, which halves $\rI$. We continue in this way, splitting along $N$ whenever possible, and along $\RR$ once it is not, always selecting the unique child containing $\mathbf x$, as long as the descendant obtained remains large. This process cannot continue infinitely many times without effect on $\rI$: after finitely many such splittings, either $\rI$ has already dropped below $\varepsilon$ (in which case we stop and are done) or we reach a transition admissible set, at which point $\rI$ equals $1/2$.
     
    From this transition point on (or from the start, if $R$ was already small), the descendant is small, and from here the splitting proceeds dyadically: at each step the admissible set has finitely many children, each with $\rI$ halved, and again exactly one of them contains $\mathbf x$. Iterating this dyadic splitting, still always choosing the child containing $\mathbf x$, produces a sequence of descendants along which $\rI$ goes to $0$. Hence, after finitely many steps, we reach a descendant $R'$ containing $\mathbf{x}$ such that $\rI'\leq\varepsilon$, as we wanted. 
\end{proof}
We are now ready to prove Theorem \ref{Thm-CZspaceproof}, which states that $S$, endowed with the left-invariant Riemannian distance $\riemdist$, the right Haar measure $\rho$ and the system of admissible sets described above $\mathcal A$, is a regular CZ space in the sense of \cite[Definition 4.1]{CMV26}.

\begin{proof}[Proof of Theorem \ref{Thm-CZspaceproof}.]
\ref{regular a} Comes from the fact that any admissible set is measurable and bounded.\\
\ref{regular b} For the exhaustion property, choose \(R_1\in\admis\) such that
\[
\overline{\damekricciball{e}{1}}\subseteq R_1,
\]
which is possible by Corollary \ref{Rmk-compactinadmissible}. Suppose that \(R_n\) has been chosen. Since \(R_n\) is bounded and \(S\) is connected and complete, the set
\[
K_n\defeq\overline{R_n}\cup\overline{\damekricciball{e}{n+1}}
\]
is compact. Corollary \ref{Rmk-compactinadmissible} therefore yields an admissible set
\(R_{n+1}\) such that \(K_n\subseteq R_{n+1}\). Thus \(R_n\subseteq R_{n+1}\) and
\(\damekricciball{e}{n}\subseteq R_n\) for every \(n\), and therefore
\[
S=\bigcup_{n\in\NN}R_n.
\]
\ref{regular c} Follows from the construction of children and descendants described above. Since every admissible set has only finitely many children, also the family \(\children^*(R)\) of all descendants of \(R\) is at most countable.\\
 \ref{regular d} Follows from Corollary \ref{cor-R contenuto in palla} and Proposition \ref{prop-discendenti piccoli}.\\
\ref{regular R} Regularity follows from Proposition \ref{prop-regularity}.\\
\ref{regular CZ} The CZ property follows from Propositions \ref{Propo-largeinclusions} and \ref{Propo-smallinclusions}.
\end{proof}

Theorem \ref{Thm-CZspaceproof} shows that \(\tonde{S,\damekriccidist,\mathcal B(S),\rho,\admis}\) is a regular CZ space. We may therefore define the associated atomic Hardy spaces \(H^{1}(S)\), and $\BMO(S)$ spaces as in \cite[Section 3]{CMV26}.
\begin{definition}\label{Def-Hardyabstract}
\begin{itemize}
    \item An $\infty$-atom on \(\tonde{S,\damekriccidist,\mathcal B(S),\rho,\admis}\) is a function of the form 
        $$
        a - \tonde{\int a \ \dd m} \frac{1}{m(R)} \indicator{ R}
        $$
         for some $R\in \admis$ and $a \in L^\infty(S)$ with $\graffe{a \neq 0 }\subseteq R$ and $\norma{a}_\infty\leq m(R)^{-1}.$
         \item We define $H^{1,\infty}(S)$ as the set of functions $f\in L^1(S)$ that can be written as $$f=\sum_{n\in\NN} \lambda_nb_n,$$ with convergence in $L^1(S)$, where the $b_n$ are $\infty$-atoms and $(\lambda_n)_{n}\in \ell^1.$ For $f\in H^{1,\infty}(S),$ we define \begin{equation*}
             \norma{f}_{H^{1,\infty}(S)}=\inf \graffe{\sum_{n\in \NN} \abs{\lambda_n} \ : \ f=\sum_{n\in \NN} \lambda_nb_n, \ b_n \ \infty\text{-atoms}, \ (\lambda_n)_n \in \ell^1}.
         \end{equation*}   
         \item We set $$L^1_{\mathrm{loc}}(S):=\graffe {f \ \text{measurable on } S \ : \ f\indicator{A} \in L^q(S) \ \forall A \in \admis} .$$ We define $BMO_1(S)$ as the set of $f\in L^1_{\mathrm{loc}}(X)$ such that 
    $$
    \norma{f}_{BMO_1}:=\sup_{ A\in \admis } \tonde{\frac{1}{m(A)}\int_A \abs{f-\tonde{\frac{1}{m(A)}\int_{ A} f \ \dd m}} \dd m} <\infty.
    $$
    The space $\BMO_1(X)$ is then defined as $BMO_1(X)$ modulo constant functions, which makes it a Banach space.
\end{itemize}
\end{definition}
\begin{remark}\label{Re-differenthardy}
   Above we have defined the Hardy space $H^{1,\infty}(S)$ and the space $\BMO_1(S)$, but we could equally have defined the Hardy spaces $H^{1,q}(S)$ and the spaces $\BMO_q(S)$ for any $q\in(1,+\infty)$. All these spaces coincide, with equivalence of norms, by \cite[Proposition 3.29]{CMV26} and \cite[Corollary 3.8]{CMV26}, respectively. We may therefore write simply $H^1(S)$ and $\BMO(S)$.
   
   The definitions of ``atom'' and of $\BMO$ spaces given above differ from \cite[Definition 3.10]{CMV26} and \cite[Definition 3.1]{CMV26}, respectively: in both cases, the average of $a$ (or of $f$) is subtracted over a parent of $A$ rather than over $A$ itself. The alternative definitions in \cite{CMV26} are needed because, on a general CZ space, there is no control on the ratio between the measure of a parent and that of its child. In the regular case, however, such control does hold, and \cite[Remark 2.28]{CMV26} and \cite[Remark 2.14]{CMV26} guarantee that the two definitions lead to the same Hardy space $H^{1}$ and the same space $\BMO$, respectively (with equivalence of norms). We shall therefore adopt the definitions that are closer to the classical case.
\end{remark}
We collect here some properties of Hardy and $\BMO$ spaces, which can be found in \cite[Proposition 3.29, Lemma 2.54]{CMV26}.
\begin{proposition}
\label{Pro-equivalence}
\label{Prop-duality}
\label{lem-density in H1}
\begin{enumerate}
On $(S, \damekriccidist, \mathcal B(S),\rho, \admis)$:
    \item if $q\in (1,\infty)$ and $p\in [q,\infty]$, then $H^1(X)\cap L^p(S)$ is dense in $H^1(S)$ and it is densely contained in $L^{q}(S);$
    \item there is an isomorphism between $\BMO(S)$ and the dual space of $H^1(S)$, with equivalence of norms, realized by the natural pairing given by integration: $g\in\BMO(S)$ acts on $f\in H^1(S)$ via $\int_S fg\,\dd\rho$ when $f$ is a finite linear combination of atoms, and this action extends uniquely by density to a bounded functional on all the functions of $H^1(S)$.
\end{enumerate}
\end{proposition}

\section{Endpoint behaviour of Riesz transforms}\label{Sec-endpointRiesz}
The expressions for operators appearing in this section are adapted from \cite{LM26}, where the multiplicative coordinate \(a\in\RR_+\) is used; here we express them in the additive coordinate \(u=\log a\), consistently with the notation adopted in the preceding sections. Under this identification, one has
\[
 a=\ee^u, \quad a\partial_a=\partial_u,\quad \frac{\dd a}{a}=\dd u.
\]

First, we recall some definitions from \cite[Section 1]{LM26}. Let \(e_1,\dots,e_{\dimV+\dimZ}\) be an orthonormal basis of \(\mathfrak n=\mathfrak v\oplus\mathfrak z\), adapted to the orthogonal decomposition. The corresponding basis \(\vfXN_1,\dots,\vfXN_{\dimV+\dimZ}\) of left-invariant vector fields on \(N\),  obtained by extending $e_1,\dots,e_{\dimV+\dimZ}$ by left translations, is given by
\begin{equation}\label{Def-leftinvariantfieldsN}
\begin{aligned}
    \vfXN_j &=\partial_{x_j} +\sum_{i=1}^{\dimZ} \frac 1 2 \quadre{x,e_j}_i \partial_{z_i}, \quad & j=1,\dots,\dimV,\\
    \vfXN_{\dimV+j}&=\partial_{z_j},
    & j=1,\dots,\dimZ.
\end{aligned}
\end{equation}
The corresponding basis of left-invariant vector fields on \(S=N\rtimes\RR\) is obtained by lifting \(X_1^\flat,\ldots,X_{\dimV+\dimZ}^\flat\) to \(S\) and adding the vector field \(X_0=\partial_u\), corresponding to the $\RR$-component. More precisely,
\begin{equation}\label{Def-leftinvariantfieldsS}
\begin{aligned}
    \vfX_0&=\partial_u,\\
    \vfX_j&=\ee^{u/2}\vfXN_j= \ee^{u/2}
    \left(\partial_{x_j} +\sum_{i=1}^{\dimZ} \frac 1 2 \quadre{x,e_j}_i \partial_{z_i}\right),
    && j=1,\dots,\dimV,\\
    \vfX_{\dimV+j}&=\ee^u\vfXN_{\dimV+j}=\ee^u\partial_{z_j},
    && j=1,\dots,\dimZ.
\end{aligned}
\end{equation}

Since left-invariant vector fields are skew-adjoint with respect to the right Haar measure, the distinguished Laplacian \(\Delta\) on \(S\) is the second-order elliptic differential operator defined by
\begin{equation}\label{Def-distinguishedlaplacian}
    \Delta= -\sum_{j=0}^{\dimV+\dimZ}\vfX_j^2,
\end{equation}
which is left-invariant and essentially self-adjoint on \(C_c^\infty(S)\subseteq L^2(S)\). The distinguished Laplacian can also be written in divergence form as
\[
    \Delta=\nabla^*\nabla,
\,\,\,
    \nabla=\bigl(\vfX_0,\dots,\vfX_{\dimV+\dimZ}\bigr),
\]
where \(\nabla\) is the Riemannian gradient and \(\nabla^*\) is its adjoint with respect to the right Haar measure \(\rho\). 

For \(j=0,\dots,\dimV+\dimZ\), the first-order Riesz transform associated with \(\Delta\) is defined by
\begin{equation}\label{Def-Riesztransform}
    \Rz_j=\vfX_j\Delta^{-1/2}.
\end{equation}

We shall also consider the corresponding adjoint operators
\[
    \Rz_j^*=\bigl(\vfX_j\Delta^{-1/2}\bigr)^*, \,\,\,j=0,\dots,\dimV+\dimZ,
\]
where the adjoint is taken in \(L^2(S,\rho)\).

It was proved in \cite[Theorem 1.1]{LM26} that the first-order Riesz transforms are of weak type \((1,1)\) and bounded on \(L^p(S)\) for every \(p\in(1,\infty)\).

The purpose of this section is to determine the endpoint mapping properties of the first-order Riesz transforms and their adjoints, as anticipated in in Theorem \ref{Thm-mainendpoint}.

\subsection{\texorpdfstring{\(\mathbf{H^1-L^1}\)}{H1-L1} boundedness of the Riesz transforms}\label{Subsec-H1L1}

The argument in \cite{LM26} was formulated in terms of an abstract CZ structure on \(S\). In particular, no concrete family of admissible sets, and hence no corresponding atomic Hardy space, was constructed there. The endpoint boundedness from \(H^1(S)\) to \(L^1(S)\) was therefore not stated explicitly.

The construction developed in the preceding sections provides such a concrete realization. 

We now apply the estimates on the kernels of $\riesz_j$'s established in \cite{LM26} to yield boundedness at this Hardy space endpoint.

\begin{proof}[Proof of (1) of Theorem \ref{Thm-mainendpoint}]

Following the notations in \cite[Theorem~4.2]{LM26}, we set 
\[
\sqrt{\pi} \Rz_j = \sqrt{\pi} \vfX_j \Delta^{-1/2} 
= \sum_{n\in\ZZ}\int_{2^n}^{2^{n+1}} t^{-1/2} \vfX_j e^{-t\Delta} \,\dd t
\eqdef \sum_{n\in\ZZ} T_n,
\]
and $k_n$ to be convolution kernel of each $T_n$. In \cite[Theorem 4.2]{LM26},  it was proved that the following estimates ((4.1-4.2) in \cite[Proposition 4.1]{LM26} or \cite[Theorem 2.3 and Remark 2.4]{MV21}) hold:
\begin{equation}\label{41}
    \int_S |k_n(\bx)|(1+2^{-n/2}|\bx|)^\varepsilon\,\dd\rho(\bx)\leq A
\end{equation}
and
\begin{equation}\label{42}
    \int_S |\nabla k_n^*(\bx)|\,\dd\rho(x)\leq B2^{-\frac n 2 },
\end{equation}
where \(f^*(\bx)=\delta(\bx)\overline{f(\bx^{-1})}\). For any fixed \(j=0,\cdots,\dimV+\dimZ\), the kernel \(K\) of the operator \(\Rz_j\) can be written as 
\[
K(\bx,\by)=\frac1{\sqrt{\pi}}\sum_{n\in\ZZ}k_n(\by^{-1}\bx)\delta(\by)\eqdef \frac1{\sqrt{\pi}}\sum_{n\in\ZZ}L_n(\bx,\by).
\]
By \eqref{41} and a change of variables, we obtain
\[
\sup_{\by\in S}\int_S|L_n(\bx,\by)|\left(\damekriccidist(\bx,\by)/2^{n/2}\right)^\varepsilon\dd\rho(\bx)\leq A.
\]
Moreover, estimate \eqref{42} yields
\[
\int_S|L_n(\bx,\by)-L_n(\bx,\bz)|\dd\rho(\bx)\lesssim B \damekriccidist(\by,\bz)2^{-n/2}.
\]
By \cite[Lemma~4.3]{CMV26}, these two estimates imply the integral H\"ormander condition \cite[Hypothesis (1.3)]{CMV26}:
\[
C_H(\Rz_j)= \sup_{A\in \admis} \esssup_{\by,\bz \in A} \int_{X\setminus A^*} \abs{K(\bx,\by)-K(\bx,\bz)} \ \dd \rho(\bx) <\infty,\,\,\,j=0,\cdots,\dimV+\dimZ.
\]
Moreover, by \cite[Remark~2.10]{CMV26}, \cite[Hypothesis (1.4)]{CMV26} is automatically satisfied since \(\tonde{S,\damekriccidist,\mathcal B(S),\rho,\admis}\) is a regular CZ space. Combining the preceding results with the \(L^p(S)\)-boundedness
of \(\Rz_j\) for \(p\in(1,\infty)\) established in \cite[Theorem~1.1]{LM26}, we conclude from \cite[Theorem~B.3]{CMV26} that \(\Rz_j\) is bounded from
\(H^1(S)\) to \(L^1(S)\) and \(\Rz_j^*\) is bounded from
\(L^\infty(S)\) to \(\mathrm{BMO}(S)\) for \(j=0,\cdots,\dimV+\dimZ\).
\end{proof}

\subsection{\texorpdfstring{\(\mathbf{H^1-L^1}\)}{H1-L1} unboundedness of the adjoint Riesz transforms}\label{Subsec-H1L1ad}
	
Having established the \(H^1(S)-L^1(S)\) endpoint result for the Riesz transforms, we now turn to their adjoints. We shall show that the corresponding endpoint boundedness fails for every \(\Rz_j^*\). The obstruction comes from the behaviour at infinity of their convolution kernels.

Corollary \ref{Rmk-compactinadmissible} ensures that every compactly supported function on \(S\) is supported in some admissible set. Consequently, any \(A\in C_c^\infty(S)\) satisfying
\[
    \int_S A\,\mathrm d\rho=0
\]
belongs to \(H^1(S)\), and can be seen as a multiple of an atom. This observation will allow us to construct explicit \(H^1(S)\)-functions in order to test the endpoint behaviour of the adjoint Riesz transforms.

We now recall from \cite{LM26} the notation for dilations, involutions, and convolution kernels on \(N\) that will be used throughout the proof. For a function \(F\) on \(N\) and \(\lambda>0\), its normalised dilation is
\begin{equation}\label{Def-Nnormaliseddilation}
    F_{(\lambda)}(n)\defeq\lambda^{-Q}F\bigl(\dil_{\lambda^{-1}}n\bigr),
\end{equation}
with dilation \(\dil\) defined in \eqref{Def-dilation}.
In coordinates, this reads
\[
    F_{(\lambda)}(x,z) =\lambda^{-Q}F(\lambda^{-1/2}x,\lambda^{-1}z).
\]
The involution on \(N\) is defined by
\begin{equation}\label{Def-Ninvolution}
    F^*(n)\defeq  \overline{F(n^{-1})}=\overline{F(-n)},
\end{equation}
since \(n^{-1}=-n\). For a function \(f\) on \(S\), we denote its \(u\)-slice by
\begin{equation*}
    f^{[u]}(n) \defeq f(n,u).
\end{equation*}
With this notation, convolution on \(S\) can be expressed in terms of convolution on \(N\) as
\begin{equation}\label{Eq-sliceconvolutionu}
(f*g)(n,u)
=\int_{\RR}\left[f^{[u']}*_N\bigl(g^{[u-u']}\bigr)_{(\ee^{u'})}\right](n)\dd u'.
\end{equation}
see \cite[(2.8)]{LM26}. The preceding notation allows us to separate the \(N\)-variable from the \(u\)-variable in convolution computations on \(S\). We now record the explicit decomposition of the convolution kernels of the adjoint Riesz transforms obtained in \cite{LM26}. The result isolates the non-integrable parts at infinity into the kernels \(K_j\), \(K_0\), and \(\widetilde K_0\), while the remaining terms are integrable at infinity. These formulas will be the starting point for the construction of compactly supported cancellative functions whose images fail to belong to \(L^1(S)\).

\begin{lemma}[{\cite[Propositions 5.4 and 5.5]{LM26}}]\label{Lem-mainkernels}
    The convolution kernels of \(\Rz_0-\Rz_0^*\) and \(\Rz_j^*\) can be expressed as 
    \begin{equation}\label{Def-Rieszkernelexpression}
    \begin{aligned}
        k_{\Rz_0-\Rz_0^*}&= -c_{\dimV,\dimZ}\bigl(K_0+\widetilde K_0\bigr)+p_0,\\
        k_{\Rz_j^*}&=-c_{\dimV,\dimZ}K_j+p_j,
        &&j=1,\dots,\dimV+\dimZ,\\
    \end{aligned}
    \end{equation}
    for some positive \(c_{\dimV,\dimZ}\) and \(p_j\) integrable at infinity for every
    \(j=0,\dots,\dimV+\dimZ\), with
    \begin{equation}\label{Def-adjointmainkernels}
    \begin{aligned}
        \widetilde K_0(n,u)&=r_0(n) \frac{ \mathbf 1_{\{u\geq 1\}}+\mathbf 1_{\{u\leq -1\}}}{u},\\
        K_0(n,u)&=\bigl((r_0)_{(\ee^u)}(n)-r_0(n)\bigr)\frac{\mathbf 1_{\{u\geq 1\}}}{u},\\
        K_j(n,u)&=r_j(n)\frac{\mathbf 1_{\{u\leq -1\}}}{u},\,\,\, j=1,\dots,\dimV+\dimZ,
    \end{aligned}
    \end{equation}
for $n\in N$ and $u\in \RR$, and, calling $n=(x,z)\in N$,
    \begin{equation}\label{Def-r0r_j}
    \begin{aligned}
        r_0(x,z)&\defeq\frac{1+|x|^2/4}{H(x,z)^{(Q+2)/2}},\\
        r_j(x,z) &\defeq\frac{\left(1+\frac{|x|^2}{4}\right)x_j-[x,e_j]\cdot z}{2H(x,z)^{(Q+2)/2}},&&j=1,\dots,\dimV,\\
        r_{\dimV+j}(x,z)&\defeq\frac{z_j}{H(x,z)^{(Q+2)/2}},&&j=1,\dots,\dimZ,
        \end{aligned}
    \end{equation}
    where    
    \begin{equation}\label{Def-H}
        H(x,z)\defeq \left(1+\frac{|x|^2}{4}\right)^2+|z|^2.
    \end{equation}
    Moreover, the functions \(r_j\) for \(j=1,\dots,\dimV+\dimZ\) can be rewritten as 
    \begin{equation}\label{Eq-rjderivativeformula}
        r_j= Q^{-1} \bigl(\vfXN_jH^{-\tfrac{\dimHomog} {2}}\bigr)^*.
    \end{equation}
\end{lemma}
The kernel decomposition in Lemma \ref{Lem-mainkernels} does not by itself identify the endpoint mapping properties of the adjoint Riesz transforms with those of the main kernels. For this purpose, we use the local-global decomposition underlying the proof of \cite[Proposition 6.3]{LM26}, together with the endpoint CZ result of \cite[Theorem B.3]{CMV26}.
\begin{lemma}\label{Corokernel} 
    Let \(T_K\) denote the convolution operator \(T_Kf=f*K.\)
    For every \(j=1,\dots,\dimV+\dimZ\), the operator \(\Rz_j^*\) is bounded from \(H^1(S)\) to \(L^1(S)\) if and only if \(T_{K_j}\) is bounded from \(H^1(S)\) to \(L^1(S)\). The analogous equivalence holds for \(\Rz_0^*\) and the convolution operator \(T_{K_0+\widetilde K_0}\).
\end{lemma}

\begin{proof}
    Let \(\zeta\in C_c^\infty(S)\) be identically one in a neighbourhood of the identity. For
    \(j=1,\dots,\dimV+\dimZ\), expression \eqref{Def-Rieszkernelexpression} gives
    \[
        k_{\Rz_j^*}=\zeta k_{\Rz_j^*}-c_{\dimV,\dimZ}K_j+c_{\dimV,\dimZ}\zeta K_j+(1-\zeta)p_j.
    \]
    Hence
    \begin{equation}\label{Eq-adjointkernelreduction}
        \Rz_j^*=T_{\zeta k_{\Rz_j^*}}-c_{\dimV,\dimZ}T_{K_j}+T_{q_j},
    \end{equation}
    where
    \[
        q_j \defeq c_{\dimV,\dimZ}\zeta K_j+(1-\zeta)p_j\in L^1(S).
    \]
    The dyadic kernels arising from the decomposition of the local operator $T_{\zeta k_{\Rz_j^*}}$ satisfy \eqref{41} and \eqref{42}, respectively, as follows from the estimates established in the proof of \cite[Theorem 4.3]{LM26}. Therefore, by the criterion established in Subsection \ref{Subsec-H1L1}, the operator $T_{\zeta k_{\Rz_j^*}}$ is bounded from $H^1(S)$ to $L^1(S)$. Moreover, since \(q_j\in L^1(S)\), the operator \(T_{q_j}\) is bounded from \(H^1(S)\) to \(L^1(S)\). The equivalence for \(j\geq1\) now follow from \eqref{Eq-adjointkernelreduction}.

    For \(j=0\), \eqref{Def-Rieszkernelexpression} similarly yields
    \[
        \Rz_0-\Rz_0^*=T_{\zeta k_{\Rz_0-\Rz_0^*}}-c_{\dimV,\dimZ}T_{K_0+\widetilde K_0}+T_{q_0},
    \]
    for some \(q_0\in L^1(S)\). As above, both \(T_{\zeta k_{\Rz_0-\Rz_0^*}}\) and \(T_{q_0}\) are bounded from \(H^1(S)\) to \(L^1(S)\). Since \(\Rz_0\) is bounded from \(H^1(S)\) to \(L^1(S)\) from Subsection \ref{Subsec-H1L1}, the required equivalence for \(\Rz_0^*\) follows.
\end{proof}
Lemma \ref{Corokernel} moves the problem from finding endpoint $H^1\to L^1$ boundedness results for the adjoint Riesz transforms to finding the same endpoint for the convolution operators associated with \(K_j\) and \(K_0+\widetilde K_0\). It therefore suffices to construct compactly supported cancellative functions whose images under these operators are not integrable. 

We shall also use the following elementary identities concerning normalised dilations, involutions, and left-invariant derivatives on \(N\), which will be used repeatedly in the convolution estimates below.
	\begin{lemma}\label{Lem-scalinginvolution}
		The following properties hold.
		
		\begin{enumerate}
			\item \label{lem-item1} For any function $F$ on \(N\) and any $\lambda>0$,
			$$
			\bigl(F_{(\lambda)}\bigr)^*=(F^*)_{(\lambda)}.
			$$
			\item \label{lem-item2} For any smooth function $F$,
			$$
			\vfXN_{j}\bigl(F_{(\lambda)}\bigr)(n)
			=\begin{cases}
				\lambda^{-1/2}\bigl(\vfXN_{j} F\bigr)_{(\lambda)}(n),
				& 1\le j\leq \dimV,\\[4pt]
				\lambda^{-1}\bigl(\vfXN_{j} F\bigr)_{(\lambda)}(n),
				& \dimV+1\leq j\le \dimV+\dimZ.
			\end{cases}
			$$
			\item \label{lem-item3}Let $Y$ be an arbitrary left-invariant vector field on $N$, defined by
			\begin{equation}\label{Def-leftinY}
			YF(n)=\left.\frac{\dd}{\dd t}\right|_{t=0}F\bigl(n\exp(tY)\bigr), \quad F \ {smooth}.
			\end{equation}
			Then for any $F,G\in C_c^\infty(N)$,
			$$
			(YF)*G^*=-\,F*(YG)^*.
			$$
		\end{enumerate}
	\end{lemma}
	
\begin{proof}
\eqref{lem-item1} Follows directly from the definitions of dilation and involution in \eqref{Def-Nnormaliseddilation} and \eqref{Def-Ninvolution}. \\
\eqref{lem-item2} Follows from the covariance of left-invariant vector fields \(\vfXN_j\) in \eqref{Def-leftinvariantfieldsN} under the automorphic dilations. More precisely, with
        \begin{equation}\label{Def-kappa}
        \kappa=
        \begin{cases}
        1/2, & 1\leq j\leq\dimV,\\
        1, & \dimV+1\leq j\leq\dimV+\dimZ,
        \end{cases}
       \end{equation}
       If \(D_{\lambda}\) denotes the pullback operator
\[
    D_{\lambda}F\defeq F\circ\dil_{\lambda},
\]
then one has
\[
    \vfXN_j\circ D_{\lambda^{-1}}=\lambda^{-\kappa}D_{\lambda^{-1}}\circ\vfXN_j,
\]
and the claimed formula follows from this identity and \eqref{Def-Nnormaliseddilation}.\\
\eqref{lem-item3} For the third identity, let \(n_1\in N\). By the definition of the
left-invariant vector field \(Y\) in \eqref{Def-leftinY} and of convolution on \(N\), we have
		$$
		\begin{aligned}
			(YF)*G^*(n_1)
			=&\int_N (YF)(n_1n_2^{-1})\,G^*(n_2)\,\dd n_2\\
			=&\int_{N}\left.\frac{\dd}{\dd t}\right|_{t=0}F(n_1n_2^{-1}\exp(tY))	\overline{G(n_2^{-1})}\,\dd n_2\\
			=&\left.\frac{\dd}{\dd t}\right|_{t=0}\int_{N}F(n_1(\exp(-tY)n_2)^{-1})	\,\overline{G(n_2^{-1})}\,\dd n_2.
		\end{aligned}
		$$
		Since \(N\) is unimodular, its Haar measure is left-invariant. Making the change of variables \(n_3=\exp(-tY)n_2\), we obtain
		$$
		\begin{aligned}
			(YF)*G^*(n_1)
			=&\left.\frac{\dd}{\dd t}\right|_{t=0}\int_{N}F(n_1n_3^{-1})\,\overline{G(n_3^{-1}\exp(-tY))}\dd n_3\\
			=&\int_{N}F(n_1n_3^{-1})\left.\frac{\dd}{\dd t}\right|_{t=0}\,\overline{G(n_3^{-1}\exp(-tY))}\dd n_3\\
			=&\int_{N}F(n_1 n_3^{-1})\,\overline{(-YG)(n_3^{-1})}\dd n_3\\
			=&-\int_{N}F(n_1n_3^{-1})\,(YG)^*(n_3)\,\dd n_3=-F*(YG)^*(n_1),
		\end{aligned}
		$$
		which completes the proof.
	\end{proof}
	We are now in a position to prove the unboundedness result \eqref{H1L1ubdd} in Theorem \ref{Thm-mainendpoint}. The identities in Lemma \ref{Lem-scalinginvolution} will be used to transfer derivatives through dilations, involutions, and convolutions in the construction of suitable cancellative test functions.
	\begin{proof}[Proof of \eqref{H1L1bdd} of Theorem \ref{Thm-mainendpoint}]
		In view of Lemma \ref{Corokernel}, for each \(j=1,\dots,\dimV+\dimZ\), it suffices to prove that the convolution operator with kernel \(K_j\) is not bounded from \(H^1(S)\) to \(L^1(S)\). By Corollary \ref{Rmk-compactinadmissible} , the support of every function in
        \(C_c^\infty(S)\) is contained in some admissible set. Hence, if \(A_j\in C_c^\infty(S)\) satisfies
       \[
       \int_S A_j\ \dd \rho=0,
        \]
        then \(A_j\in H^1(S)\). Therefore, it is enough to construct such a function \(A_j\) for which
        \[
          A_j*K_j\notin L^1(S).
        \]
		For \(j=0\), the same reduction applies to \(K_0+\widetilde{K}_0\).

        Let $\phi\in C_c^\infty(N)$ be a nonnegative test function supported in a small neighborhood of the identity.
        
        \textbf{Case I.} For \(j=1,\dots,\dimV+\dimZ\).  Define
		\begin{equation}\label{Def-bj}
		b_j=\vfXN_{j}\phi.
		\end{equation}
		Since $\vfXN_{j}$ is left-invariant vector field on $N$,
		$$
		\operatorname{supp} b_j\subseteq \operatorname{supp}\phi \quad \mathrm{and}\quad		\int_N b_j(n)\dd n=0.
		$$
		
		Let $\chi\in C_c^\infty(\RR)$, $\supp \chi \subseteq (-1,1)$, $\chi\geq 0$ in $(-1,1)$ and such that $\chi\equiv1$ on $\tonde{-\tfrac 1 2 ,\tfrac 1 2}$. Define
		$$
		A_j(\mathbf{x})=A_j(n,u)\defeq b_j(n)\,\chi(u).
		$$
		Then $A_j\in C_c^\infty(S)$ and
		$$
		\int_S A_j(\mathbf{x})\,\dd\rho(\mathbf{x})
		=\Big(\int_N b_j(n)\dd n\Big)
		\Big(\int_{-\infty}^\infty \chi(u)\,\dd u\Big)
		=0.
		$$
			Applying the identity \eqref{Eq-sliceconvolutionu}, it follows that
		\begin{equation*}
			\begin{aligned}
				A_{j}*K_{j}(n,u) = \int_{\RR} A_{j}^{[u']} *_N \left( K_{j}\right)^{[u-u']} _{(\ee^{u'})} \,\dd u'.
			\end{aligned}
		\end{equation*}
		Substituting the definitions of $A_j$ and $K_j$ in \eqref{Def-adjointmainkernels} yields
		\begin{equation}\label{Def-ajkj}
		\begin{aligned}
			(A_j*K_j)(n,u)
			&=\int_{\RR} I_{n}(u')\chi(u')J(u-u')\,\dd u'=\int_{\RR} I_{n}(u-u')\chi(u-u')J(u')\,\dd u',
		\end{aligned}
		\end{equation}
		where 
		\begin{equation*}
        I_{n}(u)=b_j*(r_{j})_{(\ee^{u})}(n),\,\,\,J(u)=\frac{\mathbf{1}_{\{u\le -1\}}}{u}.
        \end{equation*}
		For arbitrary $n=(x,z)\in N$, since $b_j$ is compactly supported, the $N$-integration takes place in a fixed compact set. According to \eqref{Eq-rjderivativeformula} and \eqref{Def-bj}, \(I_n\) can be written as
		$$
		\begin{aligned}
			I_{n}(u)=b_{j}*(r_{j})_{(\ee^{u})}(n)=Q^{-1}(\vfXN_{j}\phi)*\left( (\vfXN_{j}H^{-\tfrac{\dimHomog} {2}})^{*}\right) _{(\ee^{u})}(n). 
		\end{aligned}
		$$
		From Lemma \ref{Lem-scalinginvolution}, with \(\kappa\) from \eqref{Def-kappa}, we obtain that 
		\begin{equation}\label{Def-In}
		\begin{aligned}
			Q I_{n}(u)&=(\vfXN_{j}\phi)*\left( (\vfXN_{j}H^{-\tfrac{\dimHomog} {2}})_{(\ee^{u})}\right) ^{*}\\
			&=-\phi*\left( \vfXN_{j} (\vfXN_{j}H^{-\tfrac{\dimHomog} {2}})_{(\ee^{u})}\right)^*\\
			&=-\ee^{-\kappa u}\phi*\left( \left(\vfXN_{j} (\vfXN_{j}H^{-\tfrac{\dimHomog} {2}}) \right)_{(\ee^{u})} \right) ^{*}\\
			&=-\ee^{-\kappa u}\phi* \left(\vfXN_{j} (\vfXN_{j}H^{-\tfrac{\dimHomog} {2}}) \right)_{(\ee^{u})}^{*}.
		\end{aligned}
		\end{equation}
        
		\textbf{Subcase I.} For $j=\dimV+1,\dots,\dimV+\dimZ$, from the definitions of \(\vfXN_j\) and \(H\) in \eqref{Def-leftinvariantfieldsN} and \eqref{Def-H}, we can obtain, calling $n=(x,z) =(x_1,\dots,x_{\dimV}, z_1,\dots, z_{\dimZ})$ 
 		$$
		\begin{aligned}
			\vfXN_{j}(\vfXN_{j}H^{-\tfrac{\dimHomog} {2}})(x,z)
			=&\,\dimHomog(\dimHomog+2)z_{j-\dimV}^{2}H^{-\dimHomog/2-2}(x,z)-\dimHomog H^{-\dimHomog/2-1}(x,z)\\
			=&\,\dimHomog H^{-\dimHomog/2-2}(x,z)\left((Q+2)z_{j-\dimV}^{2}-|z|^{2}-\left(1+|x|^2/{4} \right)^{2} \right).
		\end{aligned}
		$$
		From the formula above, it is obvious that 
		$$	\vfXN_{j}(\vfXN_{j}H^{-\dimHomog/2})=\left(\vfXN_{j}(\vfXN_{j}H^{-\dimHomog/2}) \right)^*.$$
		Also, if $(x,z)$ is in the region 
		$$V_{j}:=\left\lbrace (x,z)\in N: z_{j-\dimV}^2\geq \frac{3|z|^{2}}{(\dimHomog+2)}, \,|z|> 1+|x|^2/4\right\rbrace ,$$ 
		it is easily seen that 
		$$
		\vfXN_{j} (\vfXN_{j}H^{-\tfrac \dimHomog 2 })(x,z)
		\gtrsim H^{-\tfrac{\dimHomog} {2}-2}(x,z)|z|^2.
		$$
		We now apply the scaling and involution. Since \(\operatorname{supp}\phi\) is contained in a sufficiently small neighbourhood of the identity and \(\ee^{u-u'}\in[\ee^{-1},\ee]\) whenever \(|u-u'|\leq1\), we may choose a region \(U_j\subseteq N\) of positive measure such that
\[
\dil_{\ee^{-(u-u')}}(n_1^{-1}n)\in V_j,\,\,\,\forall n\in U_j,\,\,n_1\in\operatorname{supp}\phi,\,\,|u-u'|\leq1.
\]
As a consequence, for any $n\in U_j$, any $n_1\in \supp \phi$ and \(|u-u'|\leq 1\), for \(n=(x,z)\),
		\begin{equation}\label{Est-xjbxjbge}
		\left( \vfXN_{j} (\vfXN_{j}H^{-\tfrac{\dimHomog} {2}})\right)^*_{(\ee^{u-u'})}(n_1^{-1}n)
		\gtrsim H(x,z)^{-\tfrac{\dimHomog} {2}-2}|z|^2.
		\end{equation}
		Combining this with \eqref{Def-In},
		$$
		\begin{aligned}
			-I_{n}(u-u')
			&\gtrsim 
			\,\ee^{-(u-u')}
			\int_N\phi(n_1)\,\left( \vfXN_{j} (\vfXN_{j}H^{-\tfrac{\dimHomog} {2}})\right)^*_{(\ee^{u-u'})}(n_1^{-1}n)\,\dd n_1.
		\end{aligned}
		$$
		Using the positivity of $\phi$ and the estimate \eqref{Est-xjbxjbge}, we deduce that
		$$
		-I_{n}(u-u')\gtrsim H(x,z)^{-\tfrac{\dimHomog} {2}-2}|z|^2,\,\,\,n\in U_j,\ |u-u'|\le 1.
		$$
		We now consider the convolution in the $u$–variable, from \eqref{Def-ajkj},
		\begin{equation*}
		(A_j*K_j)(x,z,u)=\int_{\RR}I_{n}(u-u')\,\chi(u-u')\,J(u')\,\dd u'.
		\end{equation*}
		Restricting the integral to 
\[
u-\frac12\leq u'\leq u+\frac12,
\]
on which \(\chi(u-u')=1\), we obtain, for \(u\leq-2\), so that $J(u')=\frac 1 {u'} $ and \(n\in U_j\),
		$$
		A_j*K_j(x,z,u)\gtrsim\int_{u-\frac 12}^{u+\frac 12}H(x,z)^{-\tfrac{\dimHomog} {2}-2}|z|^2\,|u'|^{-1}\,\dd u'
		\gtrsim H(x,z)^{-\tfrac{\dimHomog} {2}-2}|z|^2\,|u|^{-1}.
		$$
		Since $|u|^{-1}$ is not integrable on $(-\infty,-2]$, this shows that
		$A_j*K_j$ fails to be integrable on $U_j\times(-\infty,-2]$, and hence on the whole $S$.
		This completes the proof in the case $j=\dimV+1,\dots,\dimV+\dimZ$.

		\textbf{Subcase II.} For $j=1,\cdots,\dimV$, from the definitions of \(\vfXN_j\) and \(H\) in \eqref{Def-leftinvariantfieldsN} and \eqref{Def-H}, we can obtain 
		$$
		\begin{aligned}
			\vfXN_{j}H^{-\tfrac{\dimHomog} {2}}(x,z)=-\frac{\dimHomog} {2} H^{-\tfrac{\dimHomog} {2}-1}(x,z)\left( (1+|x|^{2}/4)x_j+[x,e_j]\cdot z\right). 
		\end{aligned}
		$$
		Applying $\vfXN_{j}$ again and taking the involution \eqref{Def-In},
		$$
		\begin{aligned}
			&\left(\vfXN_{j}(\vfXN_{j}H^{-\tfrac{\dimHomog} {2}})\right)^*(x,z)=\frac \dimHomog 4 H^{-\tfrac{\dimHomog} {2}-2}(x,z)\\
			&\,\cdot\left( (Q+2)\left[ (1+|x|^{2}/4)x_j-[x,e_j]\cdot z\right]^2-H(x,z)\left[x_j^2+2(1+|x|^2/4)+|[x,e_j]|^2\right]\right). 
		\end{aligned}
		$$
		We consider the region
		$$
		V_j=\bigl\{(x,z)\in N: |x|\geq 4,\ |z|\le |x|^2/64,\ x_j\ge 15|x|/16\bigr\}.
		$$
		On $V_j$, we have
		$$
		\left( \vfXN_{j}(\vfXN_{j} H^{-\tfrac{\dimHomog} {2}})\right)^* (x,z)\gtrsim H(x,z)^{-\tfrac{\dimHomog} {2}-2}|x|^6.
		$$
        Here we used the relations
        \[
        |[x,e_j]\cdot z|\le |x||z|, \,\,\,|[x,e_j]|^2\leq |x|^2-x^2_j,
        \]
        which can be deduced from \eqref{eq:Htype} and \eqref{Relation-Htype1}.
        
		  Following the same argument above as in the case where $j=\dimV+1,\cdots,\dimV+\dimZ$, we can determine a region $U_j\subseteq N$ also such that for any $n\in U_j$, any $n_1\in \supp \phi$ and \(|u-u'|\leq 1\),
		$$
		\left( \vfXN_{j} (\vfXN_{j}H^{-\tfrac{\dimHomog} {2}})\right)^*_{(\ee^{u-u'})} (n_1^{-1}n)
		\gtrsim H(x,z)^{-\tfrac{\dimHomog} {2}-2}|x|^6,
		$$ 
		which yields, for $n\in U_j$ and $| u-u'|\leq 1$,
		$$
		\begin{aligned}
			-I_{n}(u-u')
			&\gtrsim\ee^{-(u-u')/2}
			\int_N\phi(n_1)\left( \vfXN_{j} (\vfXN_{j}H^{-\tfrac{\dimHomog} {2}})\right)^*_{(\ee^{u-u'})} (n_1^{-1}n)\dd n_1\\
            &\gtrsim H(x,z)^{-\tfrac{\dimHomog} {2}-2}|x|^6.
		\end{aligned}
		$$
		By \eqref{Def-ajkj} and the support condition on \(\chi\), the integrand vanishes unless \(u-1<u'<u+1\). For \(u\leq-2\), we further restrict the integration to \(u-\frac12\leq u'\leq u+\frac12.\)
       On this interval, \(\chi(u-u')=1\) and \(u'\leq-\frac32\). Therefore, using the pointwise estimate above, we obtain, for \( n\in U_j\) and \(u\leq -2\),
       \[
        (A_j*K_j)(x,z,u)\gtrsim H(x,z)^{-\tfrac{\dimHomog} {2}-2}|x|^6 \int_{u-1/2}^{u+1/2}|u'|^{-1}\,\dd u'\gtrsim H(x,z)^{-\tfrac{\dimHomog} {2}-2}|x|^6|u|^{-1}.
       \]
		This estimate shows that $A_j*K_j$ fails to be integrable on $U_j\times(-\infty,-2]$ for $j=1,\cdots,\dimV$, which is sufficient for our purposes.
		
		\textbf{Case II.} We now consider \(j=0\). Here, take $b_0=b_{\dimV+1}$.
        
        Since $K_0$ is supported in the region where $u\geq 1$ by \eqref{Def-adjointmainkernels}, we deduce that $A_0 *K_0(x,z,u)=0$ whenever $u\leq -2$.
        
        Thus, by restricting to $u\leq -2$, similarly to \eqref{Def-ajkj}, we have 
		$$
		\begin{aligned}
			&A_0*(\widetilde{K}_0+K_0)(x,z,u)=A_0*\widetilde{K}_0(x,z,u)\\
			\eqdef&\int_{\RR} \widetilde{I}_{n}(u')\chi(u')J(u-u')\,\dd u'=\int_{\RR} \widetilde{I}_{n}(u-u')\chi(u-u')J(u')\,\dd u',
		\end{aligned}
		$$
		where \(J\) is as in \eqref{Def-In} and
		$$
		\begin{aligned}
			\widetilde{I}_{n}(u)=&(\partial_{z_1}\phi)*(r_0)_{(\ee^{u})}=(\partial_{z_1}\phi)*((r_0)_{(\ee^{u})})^*\\
			=&-\phi*(\partial_{z_1}(r_0)_{(\ee^{u})})^*=-\ee^{- u}\phi*\left( \left(\partial_{z_1}r_0\right)_{(\ee^{u}) }\right) ^*=-\ee^{- u}\phi* \left(\partial_{z_1}r_0\right)^*_{(\ee^{u}) }.
		\end{aligned}
		$$
		Recalling the definition of \(r_0\) in \eqref{Def-r0r_j} and \eqref{Def-Ninvolution},we have
		$$
		\left( \partial_{z_1}r_0\right)^* (x,z)=(Q+2)(1+|x|^2/4)H^{-\tfrac{\dimHomog} {2}-2}(x,z)z_{1}.
        $$
		Choosing the region $V_{0}:=\{(x,z)\in N: z_{1}>|z|/2\}$, we determine a region $U_0$ as before such that for any $n\in U_0$ and $n_1\in \supp \phi$,
		$$
		\left(\partial_{z_1}r_0\right)^*_{(\ee^{u-u'})}(n_1^{-1}n) \gtrsim H^{-\tfrac{\dimHomog} {2}-2}(1+|x|^2/4)|z|
		$$
		holds for $|u-u'|\le 1$. Then,
		$$
		\begin{aligned}
			-\widetilde{I}_{n}(u-u')&\gtrsim\,\ee^{-(u-u')}\int_N\phi(n_1)\left(\partial_{z_1}r_0\right)^*_{(\ee^{u-u'})}(n_1^{-1}n)
			\,\dd n_1\\
			&\gtrsim H(x,z)^{-\tfrac{\dimHomog} {2}-2}(1+|x|^2/4)|z|,\,\,\,n\in U_0,\ |u-u'|\le 1.
		\end{aligned}
		$$ 
		Restricting again to \(u-\frac12\leq u'\leq u+\frac12\), we obtain, for \(u\leq-2\) and \(n\in U_0\),
        \[
        \begin{aligned}
            (A_0*\widetilde K_0)(x,z,u)
            &\gtrsim H(x,z)^{-\tfrac{\dimHomog} {2}-2} (1+|x|^2/4)|z| \int_{u-1/2}^{u+1/2}|u'|^{-1}\,\dd u'\\
            &\gtrsim H(x,z)^{-\tfrac{\dimHomog} {2}-2}(1+|x|^2/4)|z||u|^{-1},
        \end{aligned}
        \]
		which is not integrable on $U_0\times (-\infty, -2]$.
		\end{proof}
	\begin{remark}
    We note that the method used here to establish the unboundedness of the adjoint Riesz transforms from \(H^1(S)\) to \(L^1(S)\) on Damek--Ricci groups is fairly general. In particular, it can be adapted to the group $ax+b$ without substantial changes to the argument. Therefore, this approach provides an alternative proof of \cite[Proposition 7.1]{M23}.
\end{remark}
	
\subsection{\texorpdfstring{\(\mathbf{L^\infty-\text{BMO}}\)}{Linfty-BMO} unboundedness of the Riesz transforms}
We conclude this section by deriving the dual counterpart of the unboundedness result established in Subsection \ref{Subsec-H1L1ad}. The \(H^1\)-\(\mathrm{BMO}\) duality from \cite[Proposition 2.39]{CMV26}, stated in Proposition \ref{Prop-duality} transfers the failure of \(H^1(S)\) to \(L^1(S)\) boundedness for the adjoint Riesz transforms to the failure of \(L^\infty(S)\) to \(\mathrm{BMO}(S)\) boundedness for the Riesz transforms. We first record the general duality principle that will be used.
	\begin{theorem}\label{Thm-H1BMOduality}
     Let \(T\) be a bounded operator on \(L^2(S)\), and denote by \(T^*\) its \(L^2\)-adjoint. Then \(T\) extends to a bounded operator from \(H^1(S)\) to \(L^1(S)\) if and only if \(T^*\) extends to a bounded operator from \(L^\infty(S)\) to \(\mathrm{BMO}(S)\).
\end{theorem}
	\begin{proof}
		We use the duality relations
		$$
		(H^1(S))^*=\mathrm{BMO}(S),
		\,\,\,
		(L^1(S))^*=L^\infty(S),
		$$
		where the first identity is from Proposition \ref{Prop-duality}, with the duality pairing given by
		$$
		\langle f,g\rangle=\int_{S} f(\mathbf{x})g(\mathbf{x})\,\dd \rho(\mathbf{x})
		$$ whenever $g\in L^\infty(S).$
		
		First, assume that $T:H^1(S)\to L^1(S)$ is bounded. Let $g\in L^\infty(S)$ and define a linear functional on $H^1(S)$ by
		$$
		\Lambda_g(f)=\left\langle Tf,g\right\rangle.
		$$
		By H\"older's inequality and the boundedness of $T$ from \(H^1(S)\) to \(L^1(S)\),
		$$
		|\Lambda_g(f)|
		\le \|Tf\|_{L^1(S)}\|g\|_{L^\infty(S)}
		\lesssim\|f\|_{H^1(S)}\|g\|_{L^\infty(S)}.
		$$
		Hence $\Lambda_g$ is a bounded linear functional on $H^1$. By Proposition \ref{Prop-duality}, there exists \(h_g\in\mathrm{BMO}(S)\), unique modulo constants, such that
	    $$
		\Lambda_g(f)=\left\langle  f,h_g\right\rangle 
		\,\,\,\text{for all } f\in H^1(S),
		$$
		and
		$$
		\|h_g\|_{\mathrm{BMO}(S)}\lesssim \|\Lambda_{g}\|=\sup_{\|f\|_{H^1(S)\le1}}|\Lambda_{g}(f)|\lesssim\|g\|_{L^\infty(S)}.
		$$
		For \(g\in L^\infty(S)\cap L^2(S)\), this element agrees with \(T^*g\) in the distributional sense. Hence \(T^*\) extends to a bounded operator from \(L^\infty(S)\) to \(\mathrm{BMO}(S)\).

		Conversely, assume that $T^*:L^\infty(S)\to\mathrm{BMO}(S)$ is bounded. 
		Let \(f\in H^1(S)\cap L^2(S)\). For every \(g\in L^\infty(S)\cap L^2(S)\), the \(L^2\)-adjoint identity and the \(H^1\)--\(\mathrm{BMO}\) duality give
		$$
		|\langle Tf, g\rangle| =|\langle f,T^*g\rangle|\lesssim\|f\|_{H^1(S)}\|T^*g\|_{\mathrm{BMO}(S)}\lesssim \|f\|_{H^1(S)}\|g\|_{L^{\infty(S)}}.
		$$
        Since \(T\) is bounded on \(L^2(S)\), we have \(Tf\in L^2(S)\). By the dual characterization of the \(L^1\)-norm on the \(\sigma\)-finite measure space \((S,\rho)\),
        \[
         \|Tf\|_{L^1(S)}=\sup_{\substack{g\in L^\infty(S)\cap L^2(S)\\ \|g\|_{L^\infty(S)}\leq1}}|\langle Tf,g\rangle|\lesssim\|f\|_{H^1(S)}.
        \]
        Since \(H^1(S)\cap L^2(S)\) is dense in \(H^1(S)\) by Proposition \ref{lem-density in H1}, the operator \(T\) extends uniquely to a bounded operator from \(H^1(S)\) to \(L^1(S)\).
	\end{proof}
    Applying Theorem \ref{Thm-H1BMOduality} to \(T=\Rz_j^*\), with \((\Rz_j^*)^*=\Rz_j\) and using the \(H^1(S)\) to \(L^1(S)\) unboundedness established in \eqref{H1L1ubdd} of Theorem \ref{Thm-mainendpoint} immediately yields \eqref{BMOubdd} of Theorem \ref{Thm-mainendpoint}.

\section*{Acknowledgements}
The authors are deeply grateful to Alessio Martini and Maria Vallarino for their valuable suggestions and insightful discussions during the preparation of this work.


\begin{thebibliography}{10}

\bibitem{ADY96}
J.-P. Anker, E. Damek and C. Yacoub. Spherical analysis on harmonic \(AN\) groups. \textit{Ann. Scuola Norm. Sup. Pisa Cl. Sci. (4)} \textbf{23(4)} (1996): 643--679.

\bibitem{CZ52}
A. P. Calder\'on and A. Zygmund. On the existence of certain singular integrals. \textit{Acta Math.} \textbf{88} (1952): 85--139.

\bibitem{CCM19}
V. Casarino, P. Ciatti and A. Martini. From refined estimates for spherical harmonics to a sharp multiplier theorem on the Grushin sphere. \textit{Adv. Math.} \textbf{350} (2019): 816--859.

\bibitem{C90}
M. Christ. A \(T(b)\) theorem with remarks on analytic capacity and the Cauchy integral. \textit{Colloq. Math.} \textbf{60/61(2)} (1990): 601--628.

\bibitem{C74}
R. R. Coifman. A real variable characterization of \(H^p\). \textit{Studia Math.} \textbf{51} (1974): 269--274.

\bibitem{CW71}
R. R. Coifman and G. Weiss. \textit{Analyse harmonique non-commutative sur certains espaces homog\`enes: \'Etude de certaines int\'egrales singuli\`eres}. Lecture Notes in Math. \textbf{242}, Springer, 1971.

\bibitem{CW77}
R. R. Coifman and G. Weiss. Extensions of Hardy spaces and their use in analysis. \textit{Bull. Amer. Math. Soc.} \textbf{83} (1977): 569--645.

\bibitem{CMV26}
J. M. Conde Alonso, A. Martini and M. Vallarino. Nondoubling Calder\'on--Zygmund theory across geometries. \textit{arXiv preprint} arXiv:2611.18991, 2026.

\bibitem{CDKR91}
M. Cowling, A. H. Dooley, A. Kor\'anyi and F. Ricci. $H$-type groups and Iwasawa decompositions. \textit{Adv. Math.} \textbf{87(1)} (1991): 1--41.

\bibitem{CDKR98}
M. Cowling, A. H. Dooley, A. Kor\'anyi and F. Ricci. An approach to symmetric spaces of rank one via groups of Heisenberg type. \textit{J. Geom. Anal.} \textbf{8(2)} (1998): 199--237.

\bibitem{CGHM94}
M.G. Cowling, S. Giulini, A. Hulanicki, and G. Mauceri, Spectral multipliers for the distinguished Laplacian on certain groups of exponential growth,
\textit{Studia Math.} \textbf{111} (1994), 103--121.

\bibitem{CD02}
G. Crandall and J. Dodziuk. Integral structures on $H$-type Lie algebras. \textit{J. Lie Theory} \textbf{12(1)} (2002): 69--79.

\bibitem{D87a}
E. Damek. Curvature of a semidirect extension of a Heisenberg type nilpotent group. \textit{Colloq. Math.} \textbf{53(2)} (1987): 249--253.

\bibitem{D87b}
E. Damek, The geometry of a semidirect extension of a Heisenberg type nilpotent group, Colloq. Math. {\bf 53} (1987), no. 2, 255--268; MR0924070

\bibitem{DR92a}
E. Damek and F. Ricci. A class of nonsymmetric harmonic Riemannian spaces. \textit{Bull. Amer. Math. Soc. (N.S.)} \textbf{27(1)} (1992): 139--142.

\bibitem{DR92b}
E. Damek and F. Ricci. Harmonic analysis on solvable extensions of $H$-type groups. \textit{J. Geom. Anal.} \textbf{2(3)} (1992): 213--248.

\bibitem{GS96}
G.I. Gaudry and P. Sj\"ogren, Singular integrals on Iwasawa $NA$ groups of rank $1$, \textit{J. Reine Angew. Math.} \textbf{479} (1996), 39--66.

\bibitem{HS03}
W. Hebisch and T. Steger. Multipliers and singular integrals on exponential growth groups. \textit{Math. Z.} \textbf{245(1)} (2003): 37--61.

\bibitem{HK12}
T. Hyt\"onen and A. Kairema. Systems of dyadic cubes in a doubling metric space. \textit{Colloq. Math.} \textbf{126(1)} (2012): 1--33.

\bibitem{H74}
A. Hulanicki, Subalgebra of $L^1(G)$ associated with Laplacian on a Lie group,
\textit{Coll. Math.} \textbf{31} (1974), 259--287.

\bibitem{K90}
A. G. Kaplan, Fundamental solutions for a class of hypoelliptic PDE generated by composition of quadratic forms, \textit{Trans. Amer. Math. Soc.} {\bf 258} (1980), no. 1, 147--153.

\bibitem{LM26}
J. Liu and A. Martini. Riesz transforms for the distinguished Laplacian on Damek--Ricci spaces and operator-valued multivariate spectral multipliers. \textit{J. Funct. Anal.} \textbf{292(2)} (2026): Article No. 111694. DOI: 10.1016/j.jfa.2026.111694. 

\bibitem{LVY12} 
L. Liu, M. Vallarino and D. C. Yang. Dyadic sets, maximal functions and applications on $ax+b$-groups, \textit{Math. Z.} {\bf 270} (2012), no. 1-2, 515--529.

\bibitem{M23}
A. Martini. Riesz transforms on \(ax+b\) groups. \textit{J. Geom. Anal.} \textbf{33(7)} (2023): Article No. 222.

\bibitem{MOV18}
A. Martini, A. Ottazzi and M. Vallarino. Spectral multipliers for sub-Laplacians on solvable extensions of stratified groups. \textit{J. Anal. Math.} \textbf{136(1)} (2018): 357--397.

\bibitem{MP25}
A. Martini and P. Plewa. Singular integrals on \(ax+b\) hypergroups and an operator-valued spectral multiplier theorem. \textit{Ann. Sc. Norm. Super. Pisa Cl. Sci.}, 2025. 

\bibitem{MV21}
A. Martini and M. Vallarino. Riesz transforms on solvable extensions of stratified groups. \textit{Studia Math.} \textbf{259(2)} (2021): 175--200.

\bibitem{MV10}
D. M\"uller and M. Vallarino, Wave equation and multiplier estimates on Damek--Ricci spaces, \textit{J. Four. Anal. Appl.} \textbf{16} (2010), 204--232.

\bibitem{NTV02}
F. L. Nazarov, S. Treil and A. L. Volberg, Accretive system $Tb$-theorems on nonhomogeneous spaces, \textit{Duke Math. J.} {\bf 113} (2002), no. 2, 259--312; MR1909219

\bibitem{S92}
R. S. Strichartz. Self-similarity on nilpotent Lie groups. In \textit{Geometric Analysis (Philadelphia, PA, 1991)}, Contemp. Math., vol. 140, Amer. Math. Soc., Providence, RI, 1992, pp. 123--157.

\bibitem{SV08}
P. Sj\"ogren and M. Vallarino. Boundedness from \(H^1\) to \(L^1\) of Riesz transforms on a Lie group of exponential growth. \textit{Ann. Inst. Fourier (Grenoble)} \textbf{58(4)} (2008): 1117--1151.

\bibitem{T01}
X. Tolsa. \(\mathrm{BMO}\), \(H^1\), and Calder\'on--Zygmund operators for non-doubling measures. \textit{Math. Ann.} \textbf{319} (2001): 89--149.

\bibitem{V07}
M. Vallarino. Spectral multipliers on Damek--Ricci spaces. \textit{J. Lie Theory} \textbf{17(1)} (2007): 163--189.

\bibitem{V09}
M. Vallarino. Spaces \(H^1\) and BMO on \(ax+b\)-groups. \textit{Collect. Math.} \textbf{60(3)} (2009): 277--295.




	\end{thebibliography}
\end{document}